\documentclass[final,12pt,times]{elsarticle}

\usepackage{graphicx}
\usepackage{multirow}
\usepackage{amssymb}
\usepackage{lineno}

\usepackage{mdframed}
\usepackage{listings}
\usepackage{url}
\usepackage{amsmath,amsfonts,amssymb}
\usepackage[ruled,vlined]{algorithm2e}
\usepackage{lipsum}
\usepackage{amsfonts}
\usepackage{graphicx}
\usepackage{epstopdf}
\usepackage{algorithmic}
\usepackage{caption}
\usepackage{subcaption}
\usepackage{fontawesome}
\ifpdf
  \DeclareGraphicsExtensions{.eps,.pdf,.png,.jpg}
\else
  \DeclareGraphicsExtensions{.eps}
\fi

\usepackage[a4paper,bindingoffset=0.0in,%
            left=0.8in,right=0.8in,top=1.in,bottom=1.in,%
            footskip=.4in]{geometry}

\usepackage{url}
\usepackage{color}
\usepackage{xcolor}
\usepackage{xspace}
\usepackage{numprint}
\usepackage{amsopn}
\usepackage{float}

\usepackage{multirow}
\usepackage{ctable}
\usepackage{rotating}
\usepackage{lscape}
\usepackage{array}
\usepackage{booktabs}
\usepackage{supertabular}
\usepackage{tabularx}
\usepackage{amsmath}
\usepackage{hyperref}

\newcommand{\vm}[1]{\mathbf{#1}}
\newcommand{\vb}[1]{\boldsymbol{#1}}

\newcommand\phivec{\boldsymbol{\phi}}

\newcommand\psivec{\boldsymbol{\psi}}

\usepackage[normalem]{ulem}

\definecolor{codegreen}{rgb}{0,0.6,0}
\definecolor{codegray}{rgb}{0.5,0.5,0.5}
\definecolor{codepurple}{rgb}{0.58,0,0.82}
\definecolor{backcolour}{rgb}{0.95,0.95,0.92}

\lstdefinestyle{mystyle}{
    backgroundcolor=\color{backcolour},   
    commentstyle=\color{codegreen},
    keywordstyle=\color{magenta},
    numberstyle=\tiny\color{codegray},
    stringstyle=\color{codepurple},
    basicstyle=\ttfamily\footnotesize,
    breakatwhitespace=false,         
    breaklines=true,                 
    captionpos=b,                    
    keepspaces=true,                 
    numbers=right,                    
    numbersep=5pt,                  
    showspaces=false,                
    showstringspaces=false,
    showtabs=false,                  
    tabsize=2
}

\usepackage{array}
\newcolumntype{?}{!{\vrule width 1pt}}

\journal{Elsevier}

\begin{document}

\begin{frontmatter}


\title{Neural-network learning of SPOD latent dynamics}



\author[label1]{Andrea Lario}
\ead{alario@sissa.it}
\author[label2]{Romit Maulik}
\ead{rmaulik@anl.gov}
\author[label3]{Oliver T. Schmidt}
\ead{oschmidt@eng.ucsd.edu}
\author[label1]{Gianluigi Rozza}
\ead{gianluigi.rozza@sissa.it}
\author[label4]{Gianmarco Mengaldo\corref{cor1}}
\ead{mpegim@nus.edu.sg}
\cortext[cor1]{Corresponding author}

\address[label1]{Scuola Internazionale Superiore di Studi Avanzati (SISSA), Italy}
\address[label2]{Argonne National Laboratory (ANL), USA}
\address[label3]{University of California San Diego (UCSD), USA}
\address[label4]{National University of Singapore (NUS), Singapore}

\begin{abstract}
%We aim to reconstruct the latent space dynamics of high dimensional systems using model order reduction via the spectral proper orthogonal decomposition (SPOD).
%The proposed method is based on three fundamental steps: in the first, we compress the data from a high-dimensional representation to a lower dimensional one by constructing the SPOD latent space; in the second, we build the time-dependent coefficients by projecting the realizations (also referred to as snapshots) onto the reduced SPOD basis and we learn their evolution in time with the aid of recurrent neural networks; in the third, we reconstruct the high-dimensional data from the learnt lower-dimensional representation. The proposed method is demonstrated on two different test cases, namely, a compressible jet flow, and a geophysical problem known as the Madden-Julian Oscillation. An extensive comparison between SPOD and the equivalent POD-based counterpart is provided and differences between the two approaches are highlighted. The numerical results suggest that the proposed model is able to provide low rank predictions of complex statistically stationary data and to provide insights into the evolution of phenomena characterized by specific range of frequencies. The comparison between POD and SPOD surrogate strategies highlights the need for further work on the characterization of the interplay of error between data reduction techniques and neural network forecasts.
We aim to reconstruct the latent space dynamics of high dimensional, quasi-stationary systems using model order reduction via the spectral proper orthogonal decomposition (SPOD).
The proposed method is based on three fundamental steps: in the first, once that the mean flow field has been subtracted from the realizations (also referred to as snapshots), we compress the data from a high-dimensional representation to a lower dimensional one by constructing the SPOD latent space; in the second, we build the time-dependent coefficients by projecting the snapshots containing the fluctuations onto the SPOD basis and we learn their evolution in time with the aid of recurrent neural networks; in the third, we reconstruct the high-dimensional data from the learnt lower-dimensional representation. The proposed method is demonstrated on two different test cases, namely, a compressible jet flow, and a geophysical problem known as the Madden-Julian Oscillation. An extensive comparison between SPOD and the equivalent POD-based counterpart is provided and differences between the two approaches are highlighted. The numerical results suggest that the proposed model is able to provide low rank predictions of complex statistically stationary data and to provide insights into the evolution of phenomena characterized by specific range of frequencies. The comparison between POD and SPOD surrogate strategies highlights the need for further work on the characterization of the interplay of error between data reduction techniques and neural network forecasts.
\end{abstract}

\begin{keyword}
Dynamical systems \sep Reduced order modeling \sep Neural networks \sep Deep learning 


\end{keyword}

\end{frontmatter}

% \linenumbers

%%%%%%%%%%%%%%%%%%%%%%%%%%%%%%%%%%%%%%%%%%%%%
%
\section{Introduction}\label{sec:intro}
%
%%%%%%%%%%%%%%%%%%%%%%%%%%%%%%%%%%%%%%%%%%%%%

Deep learning strategies are gaining traction as a tool for surrogate modeling (also referred to as emulation or non-intrusive reduced order modeling) of dynamical systems~\cite{hesthaven2018non,wang2020recurrent,maulik2020time,hamzi2021simple,chen2021physics, hijazi2020data, pichi2021artificial, papapicco2021neural, demo2021extended, shah2021finite}. Data associated to dynamical systems has usually a space and a time component. For instance, one can think of the evolution in time of the temperature at various locations worldwide. The various locations where the temperature is provided constitute the space component of the data, while the time snapshots constitute the time component. A common case in computational physics is when we have access to values of a given quantity at $S$ spatial locations (or grid points), for $N_{\text{t}}$ time snapshots. In this case, the dimension of the problem is $S \times N_{\text{t}}$, where the spatial component $S$ might be extremely large. The key in surrogate modeling or emulation is to identify a suitable representation of the data that reduces the high dimensionality of $S$ to a more computationally manageable number, that we denote with $S_{r}$. The new reduced space of dimension $S_{r}$, often referred to as the \textit{latent space}, is where we seek to learn the time evolution of the system. The learnt reduced space can then be used to reconstruct the original high-dimensional space by inverting the data compression. Indeed, the three key steps that constitute the foundations of surrogate modeling or emulation are: (i) data compression from high- to lower-dimensional representation, (ii) learning the time evolution in lower-dimensional space, (iii) reconstruction of the high-dimensional data from lower-dimensional representation learnt. Given the nature of surrogate modeling, there are two errors that one needs to take into account, the (a) data compression/projection error, and the (b) learning error. The trade-off is usually to find a data compression efficient enough to make the problem computationally feasible while retaining sufficient accuracy. This balance is obviously task dependent, and may vary across disciplines and specific applications. In addition, the interplay of data compression and time-evolution learning must be accounted for carefully, as the learning algorithms adopted might not be optimal for the specific latent space built. Because of this delicate trade-off and complex interplay, it is not surprising that several learning approaches to emulate latent-space dynamics have been recently proposed in the literature. These include approaches where proper orthogonal decomposition (POD) and autoencoders have been used to compress the data, while recurrent and convolutional neural networks, among others, have been used for the time-evolution learning task~\cite{mohan2018deep, gonzalez2018deep, hasegawa2020cnn, maulik2021reduced, rahman2019nonintrusive}. Many works are currently exploring non-linear techniques for dimensional reductions relying, for example, either on unsupervised learning approaches such as generative adversarial network (GAN)~\cite{dos2021grassmannian} or on diffusion maps~\cite{lee2021parametric}.

In computational physics, especially in the context of engineering and geophysical flows, one of the most popular methods for reducing the dimensionality of the problem is the POD~\cite{berkooz1993proper,maulik2021pyparsvd}. This technique seeks to find a latent space of dimension $S_{r}$ that is composed of a basis vector $\boldsymbol{\phi}_j(\mathbf{x})$ that represents the spatial component of the data, and time coefficients $a_{j}(t)$ that represent the time component
\begin{equation}
    \mathbf{q}(\mathbf{x}, t) = \sum_{j=1}^{S_{r}} a_{j}(t)\boldsymbol{\phi}_{j}(\mathbf{x}) \;\;\; t = 1,\dots,N 
\end{equation}
where $t = 1, \dots, N$ are the time snapshots available, and the coefficient constitute the latent space. One typically computes the basis vector and the coefficients on a portion of the snapshots available, e.g., $t = 1,\dots,N_{\text{train}}$, where $N_{\text{train}} < N$. This portion of the data is usually referred to as training data, as it is commonly used to train a data-driven method -- e.g., a neural network -- to reconstruct the time evolution of the latent space (i.e., the coefficients). Once the data-driven method is trained and the evolution of the coefficients learnt, one reconstructs the original high-dimensional data using the pre-computed basis vector $\boldsymbol{\phi}$. The learnt dynamics of the latent space can then be used as a non-intrusive reduced order model on e.g., unseen coefficients. 

In this paper, we construct a surrogate model for the spectral proper orthogonal decomposition (SPOD) latent space via recurrent neural networks. In particular, we use the SPOD to reduce the dimensionality of the problem, thereby identifying the SPOD latent space through singular value decomposition and oblique projection of the modes, and use recurrent neural networks to learn the time evolution of the SPOD latent space. The SPOD-based surrogate model is then compared against more popular POD-based surrogate strategies. The key difference between the two dimensionality reduction methods, POD and SPOD, consists of the space where the basis vector is sought for. On the one hand, POD extracts the basis vector (also referred to as modes) as the eigenvectors of the (spatially weighted) anomaly covariance matrix. Therefore, POD provides a space-only picture of the problem. On the other hand, SPOD~\cite{Lumley:1970,towne2017spectral,schmidt2019spectral} is based on an eigenvalue decomposition of the estimated cross-spectral density (CSD) matrix of the data, that is we can leverage the frequency space to extract coherent patterns of interest at specific wavenumbers. Indeed, at each frequency, SPOD yields a set of optimally ranked, time-harmonic and orthogonal modes, that are coherent both in space and time. Furthermore, these modes are identical to optimal response modes of a stochastically forced non-normal linear dynamical system if the forcing consists of white noise. Hence, they provide a direct bridge for dynamical system modeling in computational physics, where the generative process of the underlying data comes from a set of differential equations. This useful property, could be used, for instance, in the context of geophysical flows, where it can open several opportunities in the analysis and physical understanding of Earth system model results~\cite{farrell1996generalized,monahan2009empirical,SardeshmukhSura}. In particular, we can construct different SPOD latent spaces depending on the set of wavenumbers that we are interested in. For instance, in geophysical flows, it is common to encounter phenomena that are characterized by a compact frequency band and whose evolution in time might be of interest. This is the case of several modes of variability, such as the Madden-Julian Oscillation, or El Ni{\~n}o–Southern Oscillation. In these cases, one might be interested in a detailed analysis of the frequency spectrum, and the evolution of the latent space associated to only a portion of the frequencies provided by the SPOD, and/or to perform data prewhitening and de-noising~\cite{nekkanti2021frequency,chu2021stochastic,ghate2020broadband}.

In this paper, we show how SPOD surrogate modeling with a frequency-selective latent space learning can be achieved, and we compare it against more standard POD data reduction methods. We also show how the POD can be reconstructed from the SPOD modes, and discuss, for a neural network architecture, the trade-off between data reduction error and learning error.
\begin{figure}[H]
    \centering
    \includegraphics[width=0.99\textwidth]{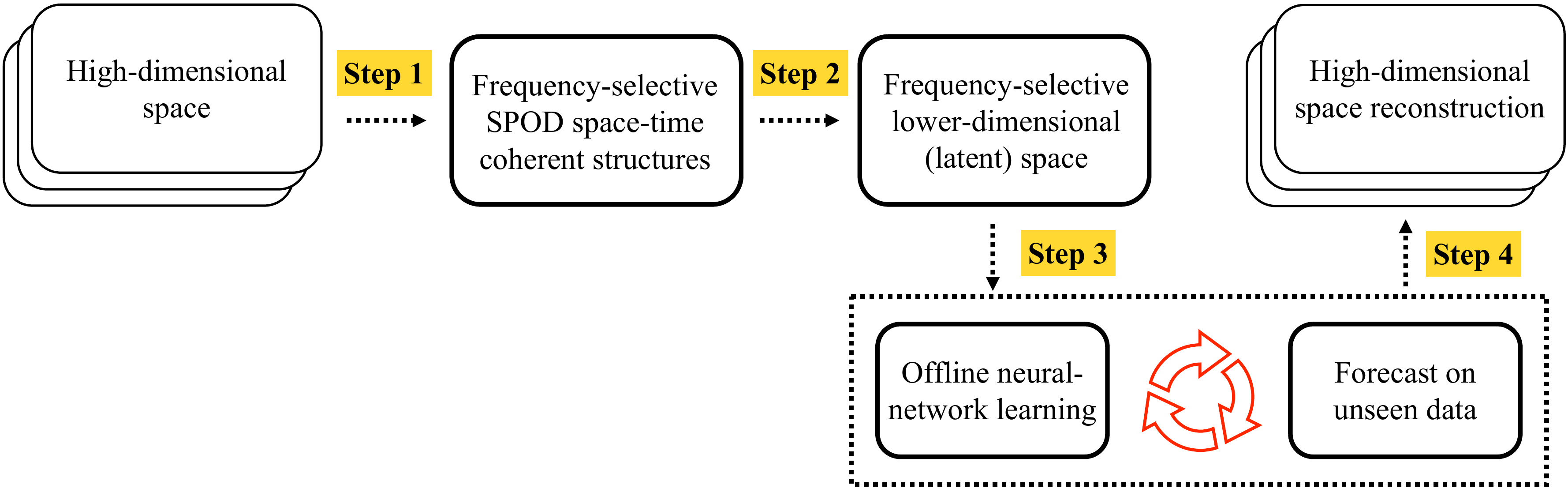}
    \caption{Conceptual workflow of SPOD latent space learning}
    \label{fig:workflow}
\end{figure}
In figure~\ref{fig:workflow}, we depict the conceptual workflow of the proposed approach for learning the latent space of SPOD, eventually selecting a frequency band of the spectrum provided by the method. The left-hand side plot of figure~\ref{fig:workflow} consists of high-dimensional data on which we apply the SPOD. From this, we obtain a set of space-time coherent SPOD structures at different periods (or frequencies), that constitute that basis to build the latent (or reduced) space that the neural network learns. We use the model that the neural network learnt to forecast future states, and we finally reconstruct high-dimensional future snapshots from the reduced space forecasts. 

The rest of the paper is organized as follows. In section~\ref{sec:methods}, we briefly introduce the SPOD method, and detail the neural-network-based SPOD emulation strategy. In section~\ref{sec:pod_vs_spod}, we compare POD and SPOD emulation, pointing out their key differences. In section~\ref{sec:results}, we present the results on a compressible jet and on a geophysical problem associated to the Madden-Julian Oscillation. We also present a direct comparison between POD and SPOD emulation. Finally, in section~\ref{sec:conclusions}, we provide a critical discussion of this work, and draw the key conclusions.

%%%%%%%%%%%%%%%%%%%%%%%%%%%%%%%%%%%%%%%%%%%%%
%
\section{SPOD emulation via non-intrusive reduced order modeling}\label{sec:methods}
%
%%%%%%%%%%%%%%%%%%%%%%%%%%%%%%%%%%%%%%%%%%%%%

This section follows the workflow presented in figure~\ref{fig:workflow}. In particular, we   introduce the SPOD method (section~\ref{sec:spod}) that allows us to reduce the dimensionality of the data and construct the SPOD latent space (section~\ref{sec:spod_latent}). These are depicted as step 1 and step 2 in figure~\ref{fig:workflow}. We then introduce the neural-network learning strategy (section~\ref{sec:learning}), depicted as step 3 in figure~\ref{fig:workflow}, along with its connections to theoretical findings. Finally, we present the reconstruction of the high-dimensional space from the SPOD latent space (section~\ref{sec:spod_reconstruction}), depicted as step 4 in figure~\ref{fig:workflow}.

\subsection{The SPOD method}\label{sec:spod}
We use the SPOD approach to identify a suitable latent space of high-dimensional data, as a first step of our dynamical system emulation strategy. The high-dimensional data consists of multidimensional  fluctuating and wide-sense stationary time-series $\mathbf{q}(\mathbf{x}, t_{i})$, from which we subtract the temporal mean $\bar{\mathbf{q}}$. The data $\mathbf{q}(\mathbf{x}, t_{i})$ is composed by $t_{i} = 1, \dots, N_{\text{t}}$ time snapshots, spatial coordinates $\mathbf{x} \in \mathbb{R}^{d}$ for a total of $S$ spatial points available, and $\mathbf{q} = q_{1}, \dots, q_{N_{\text{var}}}$ variables. 

To derive the SPOD method, we recast each time snapshot of the multidimensional data into a vector of dimension $\mathbf{q}(\mathbf{x}, t_i) = \mathbf{q}_{i} \in \mathbb{R}^{M\times N_{\text{t}}}$, where $M$ is the number of spatial points $S$ multiplied by the number of variables $N_{\text{var}}$, that is $M = S \times N_{\text{var}}$. We can then assemble the data (also referred to as snapshot) matrix
\begin{equation}
\mathbf{Q} = [\mathbf{q}_{1}, \mathbf{q}_{2}, \dots, \mathbf{q}_{N_{\text{t}}}] \in \mathbb{R}^{M\times N_{\text{t}}},
\label{eq:snapshot_matrix}
\end{equation}
that constitutes the starting point of the SPOD approach. 
The task is to identify a deterministic (set of) function(s) $\phivec(\vm{x},t)$ that best approximates the statistically stationary (zero-mean) process $\vm{q}(\vm{x},t)$ underlying our data $\mathbf{Q}$. In mathematical terms, this translates into maximizing the quantity
\begin{equation}\label{eq:SPOD_statement}
\lambda = \frac{E\{|\langle \vm{q}(\vm{x},t),\phivec(\vm{x},t) \rangle_{\vm{x},t}|^2\}}{\langle \phivec(\vm{x},t),\phivec(\vm{x},t) \rangle_{\vm{x},t}},
\end{equation}
where we assume that any realization of $\{\vm{q}(\vm{x},t)\}$ span a Hilbert space $H$ with inner product $\langle \cdot,\cdot \rangle_{\vm{x},t}$ and expectation operator $E\{\cdot\}$, here taken to be the ensemble mean. The inner product $\langle \cdot,\cdot \rangle_{\vm{x},t}$ is defined as 
\begin{equation}\label{eq:innerprod_spacetime}
\langle \vm{u},\vm{v} \rangle_{\vm{x},t} = \int_{-\infty}^{\infty} \int_{\Omega} \vm{u}^*(\vm{x},t) \vm{W} \vm{v}(\vm{x},t) \, \text{d}\Omega  \, \text{d}t, 
\end{equation}
where $\text{d} \Omega$ denotes the spatial integration, $\text{d} t$ the time integration, and $\vm{W}$ is a spatial weighting matrix. Given the assumption that any realization of $\{\vm{q}(\vm{x},t)\}$ span $H$, and thanks to Karhunen– Lo{\'e}ve (KL) theorem, we know that there exists a set of mutually orthogonal deterministic functions that forms a complete basis in the Hilbert space $H$. This can be defined as $\vm{q}(\vm{x},t) = \sum_{j=1}^{\infty} a_j\phi_{k}(\vm{x},t)$, where $\phi_{k}$ are eigenfunctions whose associated eigenvalues $\lambda_{k}$ arise from the solution of the Fredholm integral equation, that owing to the stationarity assumption on the underlying process, can be written in the frequency domain 
\begin{equation}\label{eq:SPOD}
\int_{S} \mathcal{S}(\vm{x},\vm{x}',f) \psivec(\vm{x}',f) \, \text{d}\vm{x}'  = \lambda(f) \vb{\phi}(\vm{x},f),
\end{equation}
where $\mathcal{S}(\vm{x},\vm{x}',f)$ is the Fourier transform of the two-point space-time correlation tensor $\mathcal{C}(\vm{x},\vm{x}',\tau)=E\{\vm{q}(\vm{x},t)\vm{q}^*(\vm{x}',t')\}$, also referred to as the cross-spectral density tensor
\begin{equation}\label{eq:csd}
\mathcal{S}(\vm{x},\vm{x}',f) = \int_{-\infty}^{\infty} \mathcal{C}(\vm{x},\vm{x}',\tau) e^{-\mathrm{i}2\pi f \tau} \text{d}\tau.
\end{equation}
To form the discrete analog of \eqref{eq:SPOD}, we use the Welch's method, and segment the data into $L$ (possibly overlapping) blocks  
\begin{equation}\label{eq:qmatdef}
\vm{Q}^{(\ell)} = [\vm{q}_1^{(\ell)}, \dots, \vm{q}_{N_{\text{fft}}}^{(\ell)}] \in \mathbb{R}^{M\times N_{\text{fft}}}, \;\;\; \ell = 1, \dots, L,
\end{equation}
where each block contains $N_{\text{fft}}$ snapshots, overlaps by $\ell_{\text{ovlp}}$ snapshots with the next block, and is regarded as equally representative of the whole data (i.e. it is assumed ergodic). An optimal overlap value is $50\%$ as suggested by Welch~\cite{welch1967use}, since it avoids data leakages when widowing functions are applied on the ensembles; different strategies must be adopted in same cases, for example nonoverlapping blocks minimize the variance of spectral estimations while on the other hand heavily duplicated data (i.e. high overlap value) allow to extract a higher number of blocks from the starting set of data $\mathbf{q}(\mathbf{x}, t_{i})$ for a fixed number of snapshots $N_{\text{fft}}$. The latter strategy was adopted for the test case described in \ref{sec:spod_freq_reconstruction_sensitivity} in order to extract from a relatively small data set a sufficient number of blocks to train a neural network. To obtain the discrete analog of the cross-spectral density tensor, we apply the temporal discrete Fourier transform to our data in equation~\eqref{eq:qmatdef}
\begin{equation}
\hat{\vm{Q}}^{(\ell)} = [\hat{\vm{q}}_1^{(\ell)}, \hat{\vm{q}}_2^{(\ell)}, \dots, \hat{\vm{q}}_{N_{\text{fft}}}^{(\ell)}] \in \mathbb{R}^{M \times N_{\text{fft}}}, \;\;\; \ell = 1, \dots, L.
\end{equation}
We then reorganize the Fourier-transformed data by collecting all realizations at the $k$-th frequency into a single data matrix
\begin{equation}
\hat{\vm{Q}}_k = [\hat{\vm{q}}_k^{(1)}, \hat{\vm{q}}_k^{(2)}, \dots, \hat{\vm{q}}_k^{(L)}] \in \mathbb{R}^{M \times L}, \;\;\; \text{for frequency $k$}.
\end{equation}
From this, we can construct the cross-spectral density matrix $\vm{S}_{k} = \vm{X}_k\vm{X}_k^* \in \mathbb{R}^{M \times M}$, where $\vm{X}_k = \frac{1}{\sqrt{L}}\vm{W}^{\frac{1}{2}}\hat{\vm{Q}}_k = [\vm{x}_k^{(1)},\vm{x}_k^{(2)},\dots,\vm{x}_k^{(L)}] \in \mathbb{R}^{M \times  L}$. The SPOD modes and associated energies can be computed as the eigenvectors and eigenvalues of the $\mathbf{S}_{k}$, that is $\vm{S}_{k}= \vm{U}_k \vb{\Lambda}_k \vm{U}_k^*$, where $\vb{\Lambda}_k = \mathrm{diag}(\lambda^{(1)}_{k},\lambda^{(2)}_{k}\cdots,\lambda^{(L)}_{k}) \in \mathbb{R}^{L \times L}$ and $\vm{U}_{k} = [\vm{u}^{(1)}_{k},\vm{u}^{(2)}_{k},\dots,\vm{u}^{(L)}_{k}] \in \mathbb{R}^{M \times L}$ are the diagonal matrix of eigenvalues, or modal energies, and the eigenvector matrix, respectively. In practice, we solve the more computationally amenable eigenvalue problem $\vm{X}_k^*\vm{X}_k \vm{V}_k = \vm{V}_k\vm{\Lambda}_k$ for $\vm{X}_k^*\vm{X}_k \in \mathbb{R}^{L \times L}$ instead of $\vm{X}_k\vm{X}_k^* \in \mathbb{R}^{M \times M}$, following~\cite{sirovich1987turbulence}, and obtain the SPOD modes for the $k$-th frequency from 
\begin{equation}\label{eq:discrete_spod}
\vm{\Phi}_k = \vm{W}^{-\frac{1}{2}}\vm{X}_k \vm{V}_k \vb{\Lambda}_k^{-\frac{1}{2}} = \frac{1}{\sqrt{L}} \hat{\vm{Q}}_k \vm{V}_k \vb{\Lambda}_k^{-\frac{1}{2}}.
\end{equation}
We note that if we group the frequencies together, we can rewrite the SPOD modes as follows:
\begin{equation}\label{eq:discrete_modes}
% \begin{array}{ll}
\vm{\Phi} = [\boldsymbol{\phi}_{1}, \boldsymbol{\phi}_{2}, \dots, \boldsymbol{\phi}_{N_{\text{f}}}] 
%\\[0.5em]
= [
\underbrace{\phi^{(1)}_{1}, \phi^{(2)}_{1}, \dots, \phi^{(L)}_{1}}_{\boldsymbol{\phi}_{1}}, 
\underbrace{\phi^{(1)}_{2}, \phi^{(2)}_{2}, \dots, \phi^{(L)}_{2}}_{\boldsymbol{\phi}_{2}}, \dots, 
\underbrace{\phi^{(1)}_{N_{\text{f}}}, \phi^{(2)}_{N_{\text{f}}}, \dots, \phi^{(L)}_{N_{\text{f}}}}_{\boldsymbol{\phi}_{N_{\text{f}}}}].
% \end{array}
\vspace{-0.3cm}
\end{equation}
noting that for real data the total number of frequencies is $N_{\text{f}} = \lceil \frac{N_{\text{fft}}}{2} \rceil+1$. From these modes, it is possible to construct the latent space. This is discussed next. For additional details on the SPOD method, the interested reader can refer to~\cite{towne2017spectral,schmidt2019spectral}.

\subsection{SPOD latent space}\label{sec:spod_latent}
As shown by Nekkanti and Schmidt~\cite{nekkanti2021frequency}, there are two strategies to reconstruct the SPOD latent space, namely a time-domain and a frequency-domain reconstruction. In this paper, we mainly use the time-domain reconstruction. Yet, we perform a comparison between time- and frequency-domain learning of SPOD latent space in section~\ref{sec:spod_freq_reconstruction_sensitivity}.

\subsubsection{Time-domain reconstruction}\label{sec:time-domain-reconstruction}
The SPOD latent space is constituted by a matrix of expansion coefficients $\mathbf{A}$. Following \cite{nekkanti2021frequency}, this can be constructed using a weighted oblique projection of the data onto the modal basis
\begin{equation}\label{eq:coeffs}
%\begin{array}{ll}
    \mathbf{A} = (\boldsymbol{\Phi}^{*}\mathbf{W}\boldsymbol{\Phi})^{-1}\boldsymbol{\Phi}^{*}\mathbf{W}\mathbf{Q} %\\[0.5em]
    = [
    \underbrace{a^{(1)}_{1}, a^{(2)}_{1}, \dots, a^{(L)}_{1}}_{\mathbf{a}_{1}}, 
    \underbrace{a^{(1)}_{2}, a^{(2)}_{2}, \dots, a^{(L)}_{2}}_{\mathbf{a}_{2}}, \dots, 
    \underbrace{a^{(1)}_{N_{\text{f}}}, a^{(2)}_{N_{\text{f}}}, \dots, a^{(L)}_{N_{\text{f}}}}_{\mathbf{a}_{N_{\text{f}}}}].
%\end{array}
\vspace{-0.3cm}
\end{equation}
where $\vm{A} \in \mathbb{R}^{(L\times N_{\text{f}})\times N_t}$ is the matrix containing the expansion coefficients and $\tilde{\vm{\Phi}} \in \mathbb{R}^{M \times (L\times N_{\text{f}_r})}$ is a matrix which gathers all the SPOD modes arranged by frequency as in equation~\eqref{eq:discrete_modes}. The full matrix of expansion coefficients, constructed using all modes and frequencies has dimensions $L\times N_{\text{f}}$ (frequencies). In practice, it is common to use only a portion $L_{r}$ of the total number of modes, and (eventually) a portion $N_{\text{fr}}$ of the total number of frequencies, from a frequency id $\text{flb}$ to $\text{fub}$. This recasts the original high-dimensional data into a smaller SPOD latent space of dimension $L_{r} \times N_{\text{fr}}$, defined as follows
\begin{equation}\label{eq:coeffs_r}
\begin{array}{ll}
    \mathbf{A}_r & = (\boldsymbol{\Phi}_r^{*}\mathbf{W}\boldsymbol{\Phi}_r)^{-1}\boldsymbol{\Phi}_r^{*}\mathbf{W}\mathbf{Q} = \\[0.5em]
    & = [
    \underbrace{a^{(1)}_{\text{flb}}, a^{(2)}_{\text{flb}}, \dots, a^{(L_r)}_{\text{flb}}}_{\mathbf{a}_{\text{flb}}}, 
    \underbrace{a^{(1)}_{\text{flb}+1}, a^{(2)}_{\text{flb}+1}, \dots, a^{(L_r)}_{\text{flb}+1}}_{\mathbf{a}_{\text{flb}+1}}, \dots, 
    \underbrace{a^{(1)}_{\text{fub}}, a^{(2)}_{\text{fub}}, \dots, a^{(L_r)}_{\text{fub}}}_{\mathbf{a}_{\text{fub}}}].
\end{array}
\vspace{-0.3cm}
\end{equation}
The calculation of the matrix coefficients~\eqref{eq:coeffs} can be problematic, because $\tilde{\vm{\Phi}}^* \vm{W} \tilde{\vm{\Phi}}$ is potentially ill conditioned especially when a high number of modes are taken into account. Computing its inverse numerically can therefore be difficult. In this work, we follow Nekkanti and Schmidt~\cite{nekkanti2021frequency}, and use a pseudo-inverse based on the singular value decomposition of $\tilde{\vm{\Phi}}^* \vm{W} \tilde{\vm{\Phi}}$. Additionally, we note that one cannot use an orthogonal projection because the modes at different frequencies are not orthogonal when a space inner product is considered. 

Once the reduced SPOD latent space $\mathbf{A}_r\in \mathbb{R}^{(L_r\times N_{\text{fr}}) \times N_t}$ is available, we can learn its time evolution using data-driven methods.

\subsubsection{Frequency-domain reconstruction}\label{sec:freq-domain-reconstruction}
An alternative approach for the computation of the expansion coefficients, usually referred as frequency-domain reconstruction, relies on the orthogonal projection of the Fourier-transformed data $\hat{Q}$ onto the SPOD modes:
\begin{equation}
\mathbf{\hat{A}} = \boldsymbol{\Phi}^{*}\mathbf{W}\mathbf{\hat{Q}}
\end{equation}
where $\vm{\hat{A}} \in \mathbb{R}^{(L\times N_{\text{f}})\times L}$ is the matrix containing the coefficients~\cite{nekkanti2021frequency}. 
Hence the snapshots which belong to the $l$-th block can be reconstructed by computing the corresponding Fourier transformed data starting from the modes and the expansion coefficients:
\begin{equation}
\mathbf{\hat{Q}}^{(l)} = \left[ \left(\sum_{i=1}^{L} \hat{a}_{ikl}\phi_i \right)_{k=1},  \left(\sum_{i=1}^{L} \hat{a}_{ikl}\phi_i \right)_{k=2}, ..., \left(\sum_{i=1}^{L} \hat{a}_{ikl}\phi_i \right)_{k=N_f}\right] \quad,
\end{equation}
and then by applying an Inverse Fast Fourier Transform to the vector just computed.
In order to perform properly this operation, the number of frequencies $N_f$ must be equal to the number of the snapshots contained in each block; if only some subsets of frequencies have to be retained during the snapshot reconstruction phase, one has to put equal to zero the expansion coefficients corresponding to the frequencies to be cut off.

In this paper, we propose to learn the latent space evolution via neural networks, and a long short-term memory neural network in particular. This is described next.

\subsection{Learning of dynamics via neural networks}\label{sec:learning}
We use deep learning strategies based on neural networks to learn the time evolution of the latent space constituted by the expansion coefficients constructed in equation~\ref{eq:coeffs_r}. The problem here is in the realm of autoregressive time-series forecasting, where one seek to learn the future values of a variable\footnote{In the more general case, autoregressive time-series forecasting is performed on a set of variables vs.\ just one variable} $\boldsymbol{f} = [f(t+1), f(t+2), \dots f(t+n_{\tau F})]$, provided past values of the same variable $\boldsymbol{f} = [f(t-1), f(t-2), \dots f(t-n_{\tau P})]$, where $n_{\tau F}$ can be different from $n_{\tau P}$. In our case, we seek to learn the time evolution of $\mathbf{A}_r$, that is: we want to forecast the future values $\mathbf{A}_r = [\mathbf{a}_r(t+1), \mathbf{a}_r(t+2, \dots, \mathbf{a}_r(t+n_{\tau F})]$, provided past values of these coefficients $\mathbf{A}_r = [\mathbf{a}_r(t-1), \mathbf{a}_r(t-2, \dots, \mathbf{a}_r(t-n_{\tau P})]$. To perform autoregressive time-series forecasting, using deep learning, several neural network architectures have been proposed, including recurrent neural networks, such as long short-term memory networks (LSTMs) \cite{hochreiter1997long} and gated recurrent networks \cite{chung2014empirical}, neural ordinary differential equations \cite{chen2018neural}, transformer networks \cite{jaderberg2015spatial}, and temporal convolutional networks \cite{lea2017temporal}. Since an exhaustive comparison of these architectures is beyond the scope of this article - we shall utilize the LSTM as an exemplar architecture for learning the dynamics of our observed data on the space spanned by the SPOD modes. Moreover, viable theoretical connections are relatively simple to establish for reduced-order modeling and LSTM networks.

The LSTM network is a special type of the recurrent type of neural networks. The term `recurrent' alludes to the fact that certain neuron outputs from a neural network, at each time step in the prediction, are utilized to influence predictions in future time steps. Embedding this type of inductive bias allows for the neural network to learn long and short-term correlations between instances in time from the training data set. The basic equations of the LSTM in our context for the coefficients $\mathbf{A}_r$ are given by
\begin{align}
\label{eq:LSTM_Internal}
\begin{split}
\text{input gate: }& \boldsymbol{G}_{i}=\boldsymbol{\varphi}_{S} \circ \mathcal{F}_{i}^{N_{c}}(\mathbf{A}_r), \\
\text{forget gate: }& \boldsymbol{G}_{f}=\boldsymbol{\varphi}_{S} \circ \mathcal{F}_{f}^{N_{c}}(\mathbf{A}_r), \\
\text{output gate: }& \boldsymbol{G}_{o}=\boldsymbol{\varphi}_{S} \circ \mathcal{F}_{o}^{N_{c}}(\mathbf{A}_r), \\
\text{internal state: }& \boldsymbol{s}_{t}=\boldsymbol{G}_{f} \odot \boldsymbol{s}_{t-1}+\boldsymbol{G}_{i} \odot\left(\boldsymbol{\varphi}_{T} \circ \mathcal{F}_{\mathbf{a}}^{N_{c}}(\mathbf{A})\right), \\
\text{output: }& \mathbf{h}_t = \boldsymbol{G}_{o} \circ \boldsymbol{\varphi}_{T}\left(\boldsymbol{s}_{t}\right),
\end{split}
\end{align}
where $\mathbf{A}_r \in \mathbb{R}^{(L_r\times N_{\text{fr}})\times N_{\text{t}}}$. Here, $\boldsymbol{\varphi}_{S}$ and $\boldsymbol{\varphi}_{T}$ refer to tangent-sigmoid and tangent-hyperbolic activation functions respectively (although these may be modified with a view to improving the accuracy or robustness of learning), $N_c$ is the number of hidden layer units in the specific neural network and is another modifiable parameter. Finally, $\mathcal{F}^{n}$ refers to an affine transformation and $(\odot)$ refers to a Hadamard product of two vectors. The LSTM implementation is used to advance $\mathbf{A}$ as a function of time. The interested reader may also find connections between classical implicit time-integration methods and the construction of such time-delay embedded neural architectures \cite{sherstinsky2020fundamentals}.

The most obvious connection between dynamical systems and LSTM neural networks is through the framework of time-delay embeddings for simulation of the former. A convenient starting point is the Mori-Zwanzig formalism \cite{mori1965transport,zwanzig1973nonlinear,ma2018model} which divides the variables of a dynamical system into resolved and unresolved components. Subsequently, unresolved components are `estimated' by functions that are parameterized by the resolved degrees of freedom leading to reduced system of equations for capturing the dynamics of the system. This leads to a generalized Langevin equation for the dynamical system or in other words an equation of the resolved variables with memory and noise. 

The following formulation, derived from the first step of the Mori-Zwanzig treatment, can be used to approximate the dynamics of our observed data on the reduced space
\begin{align}
\frac{d \textbf{A}_r}{d t} = e^{\mathcal{W} t} \mathcal{W} \mathbf{A}_{r},
\end{align}
where $\mathcal{W}$ is the unknown evolutionary operator of the dynamics in this space. Subsequently, we define two self-adjoint projection operators into orthogonal subspaces given by
\begin{align}
P \textbf{A}_h = \frac{(\textbf{A}_r^0, \textbf{A}_h^0)}{(\textbf{A}_r^0,\textbf{A}_r^0)} \textbf{A}_r, \quad Q = I-P,
\end{align}
with $QP=0$ and $\textbf{A}^0 = \textbf{A} (t=0)$. Therefore, $P$ may be assumed to be a projection of our full-order representation in the space spanned by SPOD modes. Subsequently, we can further expand our system as 
\begin{align}
\frac{d \textbf{A}_r}{d t} = e^{\mathcal{W} t} (Q+P) \mathcal{W} \mathbf{A}_{r},
\end{align}
which may further be decoupled to a Markov-like projection operator $\mathcal{M}$ given by 
\begin{align}
e^{\mathcal{W} t} P \mathcal{W} A_r =  \frac{(\textbf{A}_r^0, \textbf{A}_h^0)}{(\textbf{A}_r^0,\textbf{a}_r^0)} e^{\mathcal{W} t} A_r = \mathcal{M} A_r
\end{align}
and a memory operator $\mathcal{G}$ given by
\begin{align}
\label{eq:mz_memory}
e^{\mathcal{W} t} Q \mathcal{W} \mathbf{A}_r^0 = e^{Q \mathcal{W} t} Q \mathcal{W} \mathbf{A}_r^0 + \int_{0}^{t} e^{\mathcal{W}\left(t-t_{1}\right)} P \mathcal{W} e^{Q \mathcal{W} t_{1}} Q \mathcal{W} \mathbf{A}_r^0 d t_{1} = \mathcal{G} A_r,
\end{align}
from Dyson's formula \cite{evans2008statistical} with $t_1$ corresponding to the length of memory retention. The second relationship, given by Equation \ref{eq:mz_memory}, may be assumed to be a combination of memory effects and noise. The final evolution then becomes,
\begin{align}
\frac{d \textbf{A}_r}{d t} = \mathcal{G} A_r + \mathcal{M} A_r.
\end{align}

This expression can be compared to that of the internal state update within an LSTM (Equation \ref{eq:LSTM_Internal}), which we repeat below 
\begin{align}
\boldsymbol{s}_{t}=\boldsymbol{G}_{f} \odot \boldsymbol{s}_{t-1}+\boldsymbol{G}_{i} \odot\left(\boldsymbol{\varphi}_{T} \circ \mathcal{F}_{\mathbf{A}_r}^{N_{c}}(\mathbf{A}_r)\right),
\end{align}
where a linear combination of a nonlinearly transformed input vector at time $t$ with the gated result of a hidden state at a previous time $t-1$ is used to calculate the result vector at the current time. The process of carrying a state through time via gating may be assumed to be a representation of the memory integral, whereas the utilization of the current input may be assumed to be the Markovian component of the map.

\subsection{SPOD high-dimensional space reconstruction}\label{sec:spod_reconstruction}

After learning an LSTM-based surrogate model for the dynamics of the latent space, we can reconstruct predictions from this model, $\tilde{\vm{A}}_r$, by reconstructing the predicted coefficients using the SPOD modes $\vm{\Phi}_r$, i.e.,
\begin{equation}
\tilde{\vm{Q}} = \vm{\Phi}_r \tilde{\vm{A}}_r,
\label{eq:SPODrec}
\end{equation}
where one seeks $\tilde{\vm{Q}} \approx \vm{Q}$, and $\tilde{\vm{A}}_r \approx \vm{A}_r$. Indeed, from this procedure there are two types of errors that we need to consider: (1) the projection error, and the (2) learning error. The first can be easily computed by comparing the reconstruction of the training data coefficients vs.\ the true high-dimensional data prior to learning the evolution of the latent space via the neural network, i.e.,
\begin{equation}\label{eq:reconstruction_error}
    \text{err}_{\text{projection}} = \tilde{\vm{\Phi}}_r \vm{A}_r - \vm{Q}.
\end{equation}
The second can be computed using the true coefficient vs.\ the predictions of the LSTM neural-network, i.e., 
\begin{equation}\label{eq:learning_error}
    \text{err}_{\text{learning}} = \tilde{\vm{\Phi}}_r\tilde{\vm{A}}_r - \tilde{\vm{\Phi}}_r\vm{A}_r.
\end{equation}
We note that if we use only a portion of the total number of frequencies, $N_{\text{fr}}$, we do not have access to the ``true dynamics'' $\mathbf{Q}$, as this considers only the time evolution of a frequency subspace. This can be however useful in some cases. If we believe that just a portion of the frequency spectrum is relevant to the problem at hand, and the rest is, for instance, noise, hence not contributing to the dynamics of the underlying generative process, then we are evolving a physically relevant phenomenon cleaning the frequencies associated to noise. 

The total error is a combination of the two errors in equations~\eqref{eq:reconstruction_error} and \eqref{eq:learning_error}, and can be calculated as follows
\begin{equation}\label{eq:total_error}
    \text{err}_{\text{total}} = \tilde{\vm{\Phi}}_r\tilde{\vm{A}}_r - \vm{Q}.
\end{equation}
%

% In the results section (\ref{sec:results}), we will use three error metrics, $L_1$, $L_2$, and $L_{\infty}$, that are defined as follows
% %
% \begin{align}
% L_{1} = \frac{\sum_1^{} |\text{err}|}{(n_{points} n_{snap})} \hspace{1cm}
% L_{2} = \frac{\sum(|\text{err}|)^2)}{(n_{points} n_{snap})} \hspace{1cm}
% L_{\infty} = \frac{\max(|\text{err}|)}{n_{snaps}}
% \end{align}

\section{On the connections between POD and SPOD emulation}\label{sec:pod_vs_spod}
In this section, we introduce the POD method, outline how POD modes relate to SPOD modes, and detail how POD latent space differs from the SPOD one, thereby resulting in a learning task that has a different characteristics and degree of difficulty.

\subsection{The POD method}\label{subsec:pod}
Similarly to SPOD, in the POD method we seek to identify a basis $\phivec_j^{(\text{POD})}$ that best approximates the process $\mathbf{q}(\mathbf{x},t)$ underlying our data $\mathbf{Q}$, defined in equation~\eqref{eq:snapshot_matrix}. In contrast to SPOD, for POD we seek to maximize 
\begin{equation}
\lambda^{(\text{POD})} = \frac{E\{|\langle \vm{q}(\vm{x},t),\phivec^{(\text{POD})}(\vm{x},t) \rangle_{\vm{x}}|^2\}}{\langle \phivec^{(\text{POD})}(\vm{x},t),\phivec^{(\text{POD})}(\vm{x},t) \rangle_{\vm{x}}},
\end{equation}
over $\phivec^{(\text{POD})}(\vm{x},t) \in H$, where $H$ is a Hilbert space. The analogy with the SPOD modes described in the previous section is evident, but in this case the inner product $\langle \cdot, \cdot \rangle_x$ is defined only in space and therefore the two-point correlation is written as $\mathcal{C}(\vm{x},\vm{x}') = E\{|\langle \vm{q}(\vm{x},t),\phivec^{(\text{POD})}(\vm{x},t) \rangle_{\vm{x}}|^2\}$, being the dependence on time dropped. The function $\phivec(\vm{x},t) $ which maximizes the quantity $\lambda^{(\text{POD})}$ can be found by solving the eigenvalue problem: $\langle \mathcal{C}(\vm{x},\vm{x}') , \phivec^{(\text{POD})}(x') \rangle = \lambda^{(\text{POD})}\phivec^{(\text{POD})}(x)$. In terms of discrete data $\mathbf{Q}$, the eigenvalue problem reads as $\mathbf{Q}^{*}\mathbf{Q} \mathbf{V}^{(\text{POD})} = \mathbf{V}^{(\text{POD})}\boldsymbol{\Lambda}^{\text{(POD)}}$, and the POD modes can be constructed as $\boldsymbol{\Phi}^{(\text{POD})} = \mathbf{Q}\mathbf{V}^{(\text{POD})}$. Hence POD modes identify an orthonormal basis which is coherent in space and optimal in sense that it minimizes the error in terms of energy between the original snapshots~\eqref{eq:snapshot_matrix} and their orthogonal projection.

\subsection{Relation between POD and SPOD modes}
POD modes can be seen as a superimposition of space-time coherent structures, and can be directly related to SPOD modes. In order to show this, we start from the equivalence between the space correlation tensor and the space-time correlation tensor reported in Equation~\ref{eq:csd}. For a null translation in time, the latter reads
\begin{equation}\label{eq:pod_vs_spod}
\mathcal{C}(\vm{x},\vm{x}', 0) = \int_{-\infty}^{\infty}\mathcal{S}(\vm{x},\vm{x}',f)\text{d}f. 
\end{equation}
This relation shows the equivalence between a space correlation and a zero-time-lag space correlation. As described in \cite{towne2017spectral}, the cross-spectral density tensor has the following diagonal representation:
\begin{equation}
\mathcal{S}(\vm{x},\vm{x}', f) = \sum^{\infty}_{k=1} \lambda_j(f)\psivec_k(\vm{x},f) \psivec_k^*(\vm{x}',f) \;.
\label{eq:csd_diag}
\end{equation}
Putting together Equation \ref{eq:pod_vs_spod} and \ref{eq:csd_diag}, and applying the operator $\langle \cdot, \phivec(\vm{x}) \rangle$ to both sides one eventually obtains:
\begin{equation}
\lambda^{(\text{POD})} \phivec_j^{(\text{POD})}(x) = \int_{-\infty}^{\infty}  \sum^{\infty}_{k=1} \lambda_k(f)c_{jk}(f) \psivec_k(\vm{x},f) df
\end{equation}
being $c_{jk} = \langle \phivec_k(\vm{x},f), \phivec^{(POD)}_j(\vm{x})$.
This last equation shows how the POD modes can be seen as a summation of SPOD modes at many frequencies; in particular the POD modes are mixing together space-coherent dynamics which are characterized by different frequencies.

\subsection{Comparison between POD and SPOD latent space}
Using the POD modes $\boldsymbol{\Phi}^{(\text{POD})}$ defined in section~\ref{subsec:pod}, it is possible to compute the latent space for POD, that is 
\begin{equation}\label{eq:pod_coeffs_r}
\mathbf{A}_r^{(\text{POD})} = \boldsymbol{\Phi}^{(\text{POD}), *}\mathbf{Q} = [a^{(\text{POD}, 1)}, a^{(\text{POD}, 2)}, \dots, a^{(\text{POD}, L_r)}],
\end{equation}
where $\mathbf{A}_r^{(\text{POD})} \in \mathbb{R}^{L_{r} \times N_{\text{t}}}$, with $L_{r}$ being the number of modes retained to compute the coefficients (also referred to as truncation of the POD basis). If we compare equation~\eqref{eq:pod_coeffs_r} with equation~\eqref{eq:coeffs_r}, we note that, unlike POD, the SPOD latent space is richer, because it contains the frequency content of the data, that is its dimension $\mathbf{A}_r^{(\text{POD})} \in \mathbb{R}^{L_{r} \times N_{\text{t}}} \ne \mathbf{A}_r \in \mathbb{R}^{(L_r\times N_{\text{fr}})\times N_{\text{t}}}$. It is exactly this aspect that makes SPOD and POD latent spaces different. The former contains frequency information regarding the data, and is able to capture the evolution of coherent structures for finite frequency bands. The latter contains an overall representation of the data, without the frequency granularity permitted by SPOD. This difference leads to conceptually different emulation strategies, where for SPOD we can control the frequency information we want to retain in the data. This feature allows (1) filtering out unwanted noise, or (2) focusing on specific phenomena arising for given frequency contents. Obviously, the coefficients produced by POD and SPOD have a different time behaviours. Because of this, the neural-network learning task may be more difficult for one with respect to the other data reduction method. In section~\ref{sec:results}, we will show how the latent space (first two coefficients) for POD and SPOD (figure~\ref{fig:PODvsSPODcoeff}), and it will be possible to appreciate how the SPOD coefficients are less stationary if compared to the POD ones, thereby rendering the SPOD emulation task more difficult.

%%%%%%%%%%%%%%%%%%%%%%%%%%%%%%%%%%%%%%%%%%%%%
%
\section{Results}\label{sec:results}
%
%%%%%%%%%%%%%%%%%%%%%%%%%%%%%%%%%%%%%%%%%%%%%

\subsection{Data layout, neural network configuration and code availability}\label{sec:neuralNetworkConfig}
The emulation strategy proposed in this paper aims to predict future values of the SPOD coefficients $\tilde{\mathbf{A}}_{r,\text{future}}: [\mathbf{A}_{r}(t+1), \mathbf{A}_{r}(t+2), \dots, \mathbf{A}_{r}(t+n_{\tau F})]$, provided their past values $\mathbf{A}_{r,\text{past}}: [\mathbf{A}_{r}(t-1), \mathbf{A}_{r}(t-2), \dots, \mathbf{A}_{r}(t-n_{\tau P})]$. From these, we can then compute the future data snapshots $\mathbf{Q}: [\mathbf{Q}(t+1), \mathbf{Q}(t+2), \dots, \mathbf{Q}(t+n_{\tau F})]$, simply by projecting the predicted coefficients onto the pre-computed modes $\boldsymbol{\Phi}$, that is
\begin{equation}\label{eq:future_projection}
    \mathbf{Q}(t+i) = \boldsymbol{\Phi}_r\mathbf{A}_r(t+i).
\end{equation}
The coefficients $\mathbf{A}_r$ are learnt using a recurrent neural network architecture, namely the long-short term memory (LSTM) network. The coefficients learnt via LSTM are denoted by $\tilde{\mathbf{A}}_r$, as opposed to the true coefficients denoted by $\mathbf{A}_r$. The task is to learn the (generally nonlinear) map $\tilde{\mathbf{A}}_{r,\text{future}} \longrightarrow \mathcal{N}_{\text{LSTM}}(\mathbf{A}_{r,\text{past}})$ where one seek to minimize the difference between the true $\mathbf{A}_{r,\text{future}}$ and the predicted $\tilde{\mathbf{A}}_{r,\text{future}}$, under a given error norm ($L_{2}$ in our specific case). The LSTM neural network $\mathcal{N}_{\text{LSTM}}$ adopted for all the results presented in this section shared the same overall structure, and consisted of a sequential single LSTM layer, with hyperbolic tangent activation function, and dropout of 0.15. 

\begin{figure}
    \centering
    \includegraphics[width=0.9\textwidth]{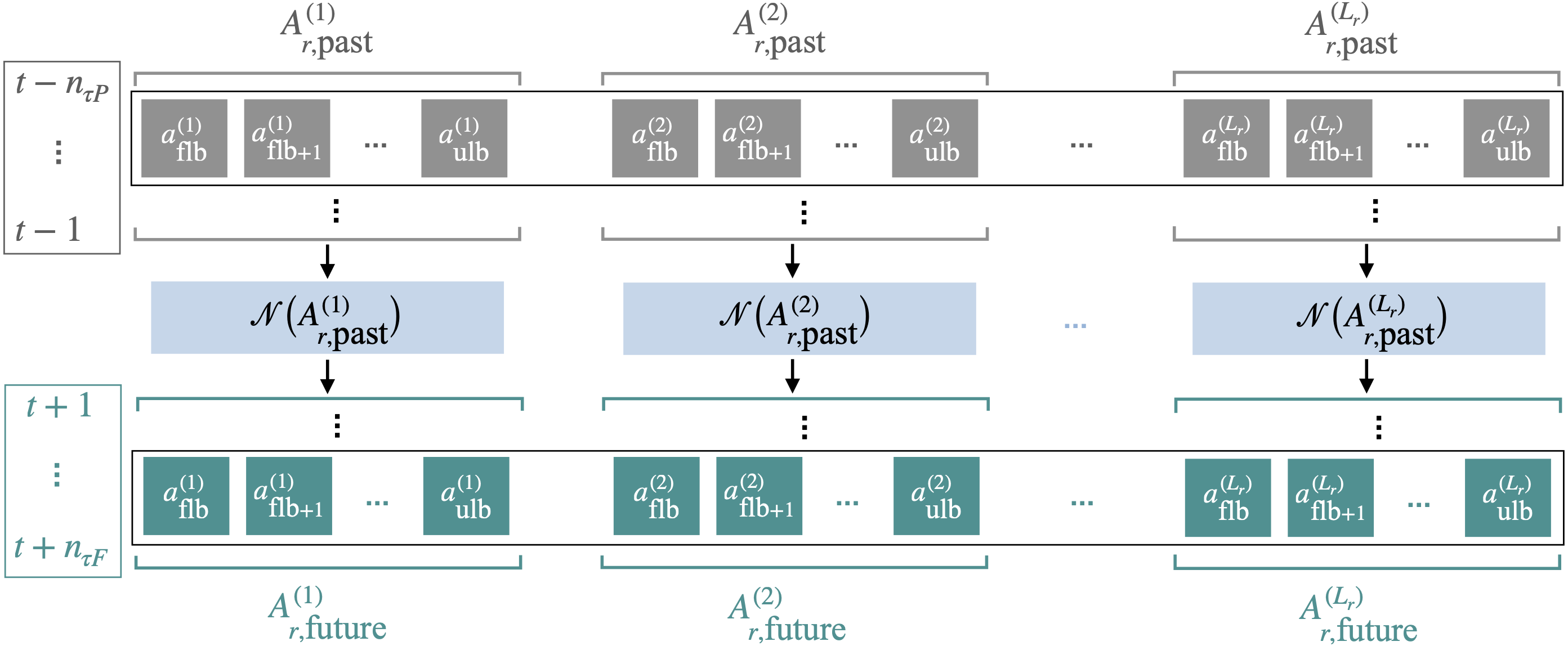}
    \caption{Learning configuration adopted for the SPOD coefficients. The grey boxes at the top represent the past values of the coefficients. The green boxes at the bottom represent the future values of the coefficients. The light blue boxes are the neural networks $\mathcal{N}$ mapping each SPOD coefficient (for all frequencies) from past to future values.}
    \label{fig:neural_nets_spod}
\end{figure}

For the learning task at hand, we divided the data available into training, validation and test sets. Consequently, the neural network was fed with
\begin{description}
    \item [a training input] containing $\mathbf{A}^{(\text{train})}_{r,\text{past}}$ of dimension $[n^{(\text{train})}_\text{{samples}}$, $n_{\tau P}$, $n_{\text{features}}]$, with $n_{\text{samples}}$ being the number of samples used for the training, $n_{\tau P}$ the number of previous time snapshots used for the future prediction, and $n_{\text{features}} = M \times N_{\text{freq}}$ the total number of features, 
    \item [a training output] containing $\mathbf{A}^{(\text{train})}_{r,\text{future}}$ of dimensions  $[n^{(\text{train})}_{\text{samples}}, n_{\tau F}\times n_{\text{features}}]$, where $n_{\tau F}$ is the number of future steps one wants to predict
    \item [a testing input] containing $\mathbf{A}^{(\text{test})}_{r,\text{past}}$ of dimension $[n^{(\text{test})}_{\text{samples}}, n_{\tau P}, n_{\text{features}}] $\
    \item [a testing output] containing $\mathbf{A}^{(\text{test})}_{r,\text{future}}$ of dimension $[n^{(\text{test})}_{\text{samples}}, n_{\tau F}\times n_{\text{features}}]$.
\end{description}
The coefficients used for testing, $\mathbf{A}^{(\text{test})}_{r}$, were computed by projecting the last $n^{(\text{test})}_{\tau P}$ data realizations on the same reduced basis (i.e., SPOD modes) obtained from the training snapshots, that is:
\begin{equation}
    \mathbf{A}^{(\text{test})}_{r} = \boldsymbol{\Phi}^{(\text{train}),*}_r\mathbf{Q}^{(\text{test})},
\end{equation}
such that information leakages which could inflate the prediction are prevented. We note that there is no validation set, as no hyperparameter optimization was performed. However, we used the test set to evaluate under- and over-fitting issues, by computing the corresponding loss function during training. All the coefficients $\mathbf{A}_{r}$ were scaled between -0.1 and 0.1 in a feature by feature fashion to avoid output saturation. To speed up neural-network training, we subdivided the coefficient matrix $\mathbf{A}_r$ into $L_{r}$ distinct sub-matrices, and distributed the computation across $L_r$ identical LSTM neural networks, where the number of features for each of them was equal to the number of frequencies used, as depicted in figure~\ref{fig:neural_nets_spod}. We note that this neural network configuration also allows for coefficient-specific hyperparameter optimization, that might be beneficial for complex learning tasks. 

We remark that we used the time-domain reconstruction detailed in section~\ref{sec:time-domain-reconstruction} for all results, except for section~\ref{sec:spod_freq_reconstruction_sensitivity}, where we compare the frequency-domain reconstruction (detailed in section~\ref{sec:freq-domain-reconstruction} vs.\ the time-domain one.

The results presented in the following sections were obtained using the open-source \textit{PySPOD} \cite{mengaldo2021pyspod}, where we implemented the neural-network emulation strategy proposed in this paper, using Keras (\url{https://keras.io}). The interested reader can refer to: \url{https://github.com/mengaldo/PySPOD}, and in particular to the Git branch \textit{papers/spod-emulation}.

\subsection{Engineering problems: a compressible jet}\label{sec:comprJet}
The first test case considered is a subsonic turbulent jet at Mach number $0.4$ and Reynolds number of $4.5\cdot10^6$. The data, retrieved from \url{https://github.com/SpectralPOD/spod_matlab}, consists of 1000 realizations (or snapshots) of a pressure field sampled every 0.2 acoustic time units and corresponds to a coarsened version of a large-eddy simulation dataset computed by Bres et al.~\cite{bres2019}. This coarsened dataset interpolates the original data onto a significantly smaller and uniformly spaced cylindrical grid of $22 \times 88$ points. The 1000 snapshots representing the pressure field at different time instants were split into a training and a testing set. The training set is composed of $n^{(\text{train})}_{\text{samples}}=800$ snapshots, while the testing set is composed of $n^{(\text{test})}_{\text{samples}}=200$. To compute the SPOD modes, the snapshots of the training set were subdivided into $L = 24$ statistically independent blocks, each of them containing 64 snapshots with a $\ell_{\text{ovlp}} = 50\%$ overlap. No weighting was adopted for the inner product, thus resulting in the weight matrix being equal to the identity matrix $\mathbf{W} = \mathbf{I}$. This choice of parameters provides up to 24 modes and 33 frequencies. The SPOD modes were generated using the methods illustrated in section \ref{sec:methods}, starting from the snapshots containing the pressure fluctuations computed by subtracting the time mean of the pressure field from each realization. In the following two subsections, we show the comparison between POD and SPOD emulation, and the sensitivity with respect to the number of modes and the number of future values predicted $n_{\tau F}$ -- also referred to as forecast horizon -- (section~\ref{sec:pod_vs_spod_sensitivity}), as well as the sensitivity of SPOD emulation to the frequencies retained (section~\ref{sec:spod_freq_sensitivity}).

\subsubsection{POD vs.\ SPOD emulation: sensitivity to number of modes and forecast horizon}\label{sec:pod_vs_spod_sensitivity}
\begin{figure}[H]
\centering
\includegraphics[width=0.6\textwidth]{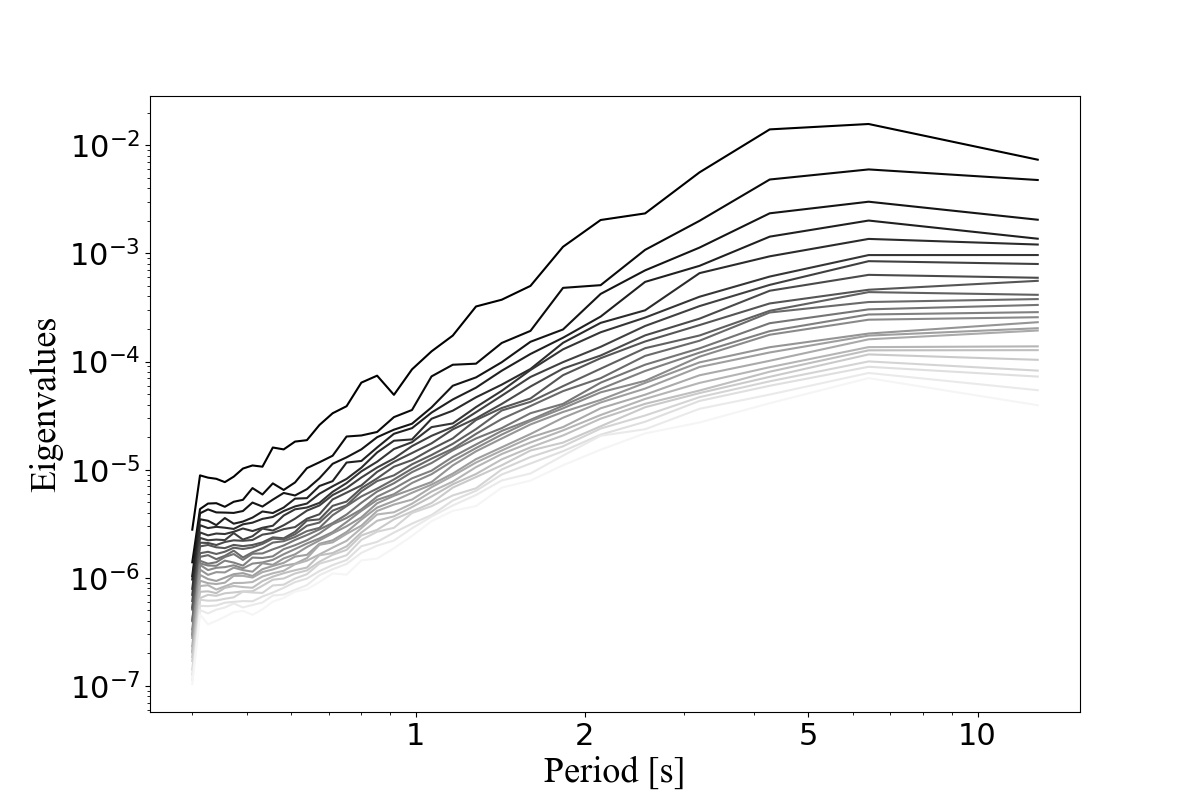}
\caption{SPOD eigenvalues plotted against period: the first few modes are able to capture most of the energy of the original system.}
\label{fig:eigVSperiod}
\end{figure}

In this section, we compare POD and SPOD emulation for different number of modes retained and different forecast horizons. We retain all 33 frequencies (equivalently periods) that the SPOD provides for each eigenvalue (and associated mode) for this test case, as depicted in figure~\ref{fig:eigVSperiod}. In contrast, the POD yields a set of eigenvalues and associated modes, that do not contain this frequency dimension. From figure~\ref{fig:eigVSperiod}, it is possible to see how the first few eigenvalues (thus associated SPOD modes) are able to capture the majority of the energy contained in the original snapshots. Therefore, as expected, we can reconstruct the dynamics of the jet pressure field by superimposing the first few leading modes multiplied by the associated coefficients computed using the oblique projection introduced in section~\ref{sec:spod_latent}. This is true for both POD and SPOD.
\begin{figure}[H]
    \centering
    \begin{subfigure}{.48\textwidth}\label{fig:recVSorig_Proj_pod}
    \includegraphics[width=\textwidth]{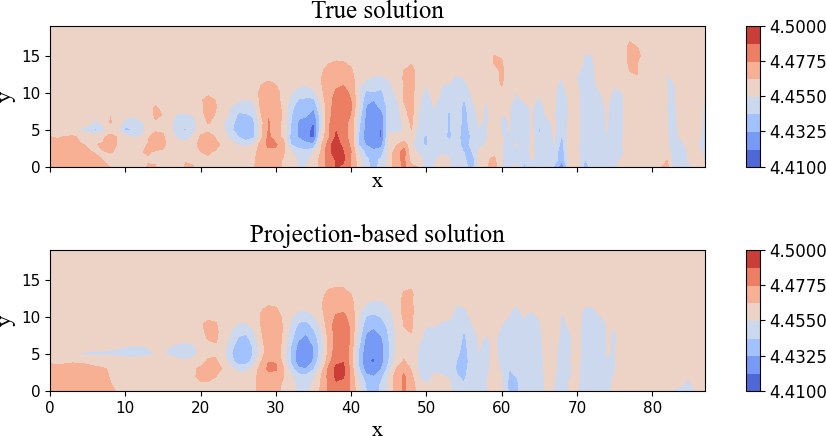}
    \caption{POD}
    \end{subfigure}
    \hspace{0.5cm}%
    \begin{subfigure}{.48\textwidth}\label{fig:recVSorig_Proj_spod}
    \includegraphics[width=\textwidth]{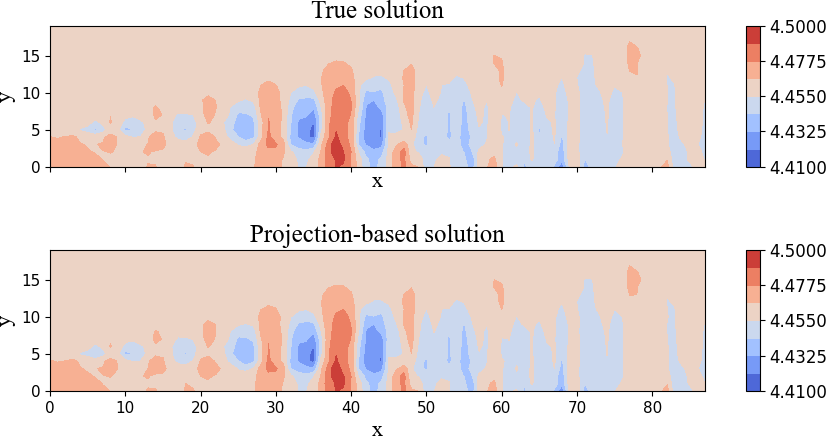}
    \caption{SPOD}
    \end{subfigure}
    \caption{True pressure field vs.\ POD reconstructed pressure field (left plot) and SPOD reconstructed pressure field (right plot). The reconstruction used ten modes for both cases (albeit the number of modes for SPOD is larger, and equal to $c = 10 (\text{modes}) \times 33 (\text{frequencies})$ as we use all 33 frequencies to reconstruct the pressure field). Both figures show time snapshot 500.}
\label{fig:recVSorig_Proj}
\end{figure}
An example of pressure field reconstruction for both methods is depicted in figure~\ref{fig:recVSorig_Proj}. Here, we show a qualitative comparison between the 500-th snapshot of the pressure field and the corresponding reconstructed solution, in which the pressure fluctuations were obtained by superimposing the first ten POD modes (left plot), and the first ten SPOD modes (right plot). 

\begin{figure}[H]
    \centering
    \begin{subfigure}{.49\textwidth}\label{fig:L2vsModes_pod}
    \includegraphics[width=\linewidth]{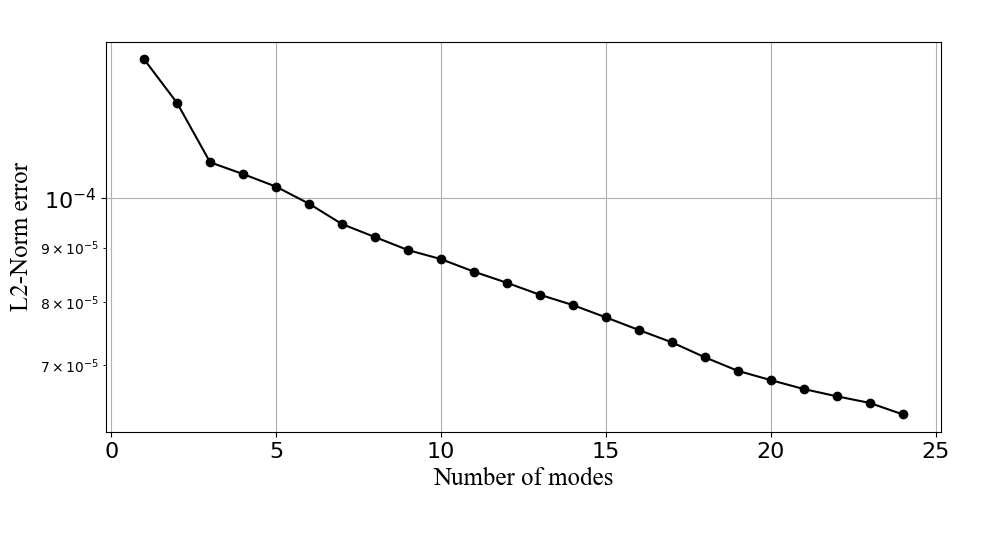}
    \caption{POD}
    \end{subfigure}
    \begin{subfigure}{.49\textwidth}\label{fig:L2vsModes_spod}
    \includegraphics[width=\linewidth]{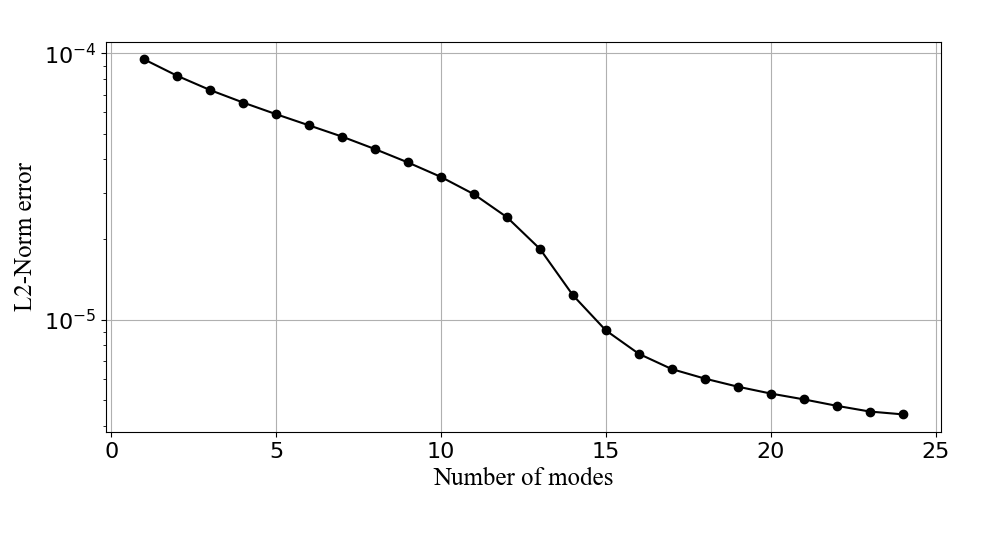}
    \caption{SPOD}
    \end{subfigure}
    \caption{Decay of the L2-norm error as a function of the number of POD (left plot) and SPOD (right plot) modes used for reconstructing the pressure field\label{fig:L2vsModes}}
\end{figure}

\begin{figure}[H]
\centering
\begin{subfigure}{.49\textwidth}
\includegraphics[width=\linewidth]{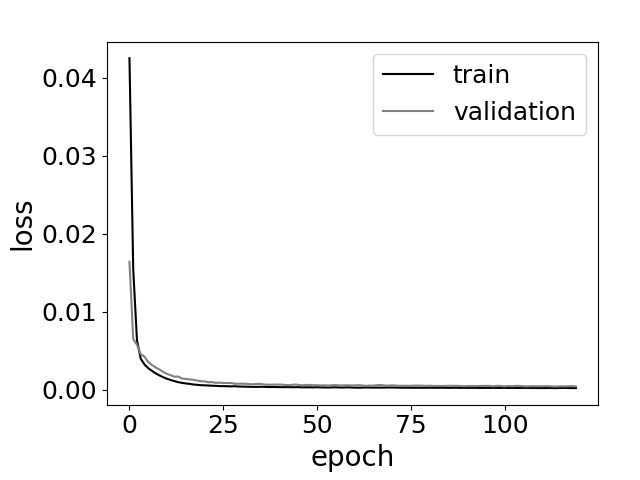}
\caption{Loss POD}
\end{subfigure}
\begin{subfigure}{.49\textwidth}
\includegraphics[width=\linewidth]{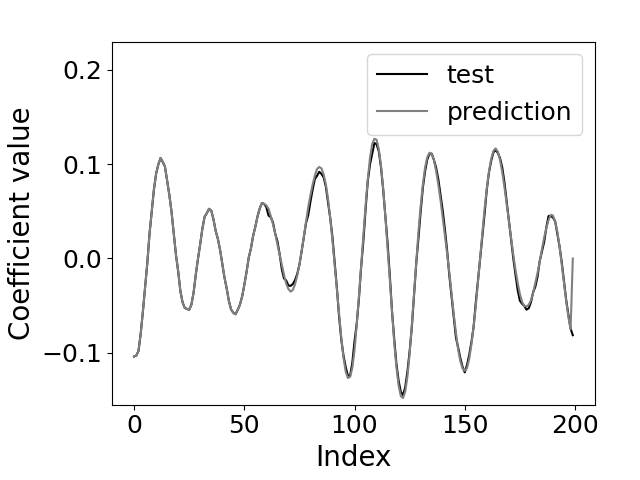}
\caption{Prediction POD}
\end{subfigure}
\begin{subfigure}{.49\textwidth}
\includegraphics[width=\linewidth]{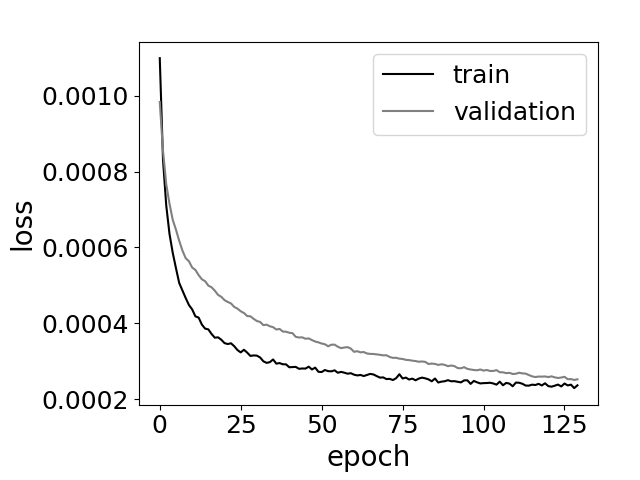}
\caption{Loss SPOD}
\end{subfigure}
\begin{subfigure}{.49\textwidth}
\includegraphics[width=\linewidth]{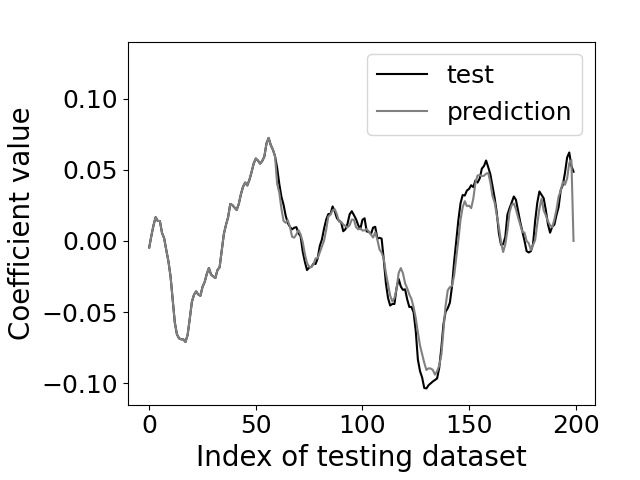}
\caption{Prediction SPOD}
\end{subfigure}
\caption{Loss histories of the training phase (right figures) and comparison between the first time coefficient and the corresponding neural network prediction (left figures). The two top figures correspond to POD, while the bottom two to the SPOD.}
\label{fig:SPODcoeff}
\end{figure}

\begin{figure}[H]
\centering
\begin{subfigure}{.49\textwidth}\label{fig:PODvsSPODcoeff_pod}
\includegraphics[width=\linewidth]{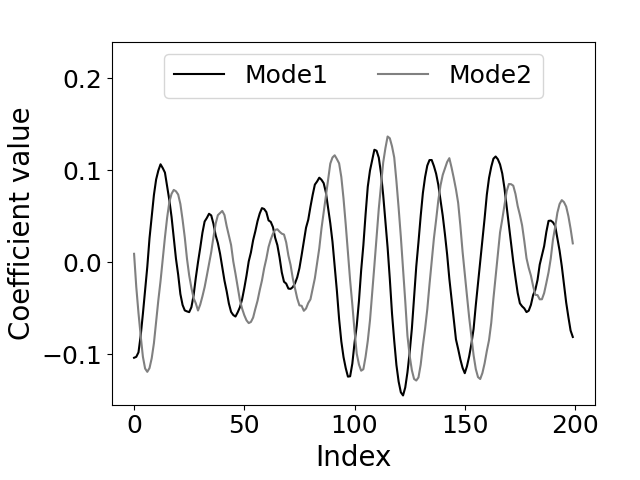}
\caption{POD}
\end{subfigure}
\begin{subfigure}{.49\textwidth}\label{fig:PODvsSPODcoeff_spod}
\includegraphics[width=\linewidth]{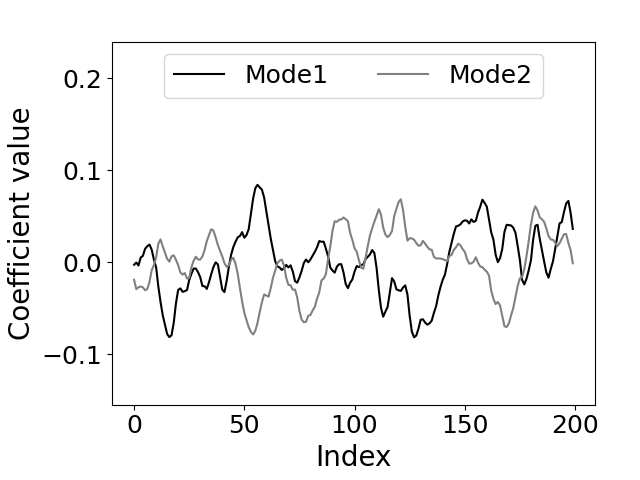}
\caption{SPOD}
\end{subfigure}
\caption{First two POD (on the left) for the same frequency and first two SPOD modes (on the right).}
\label{fig:PODvsSPODcoeff}
\end{figure}

\begin{figure}[H]
\centering
\begin{subfigure}{.48\textwidth}\label{fig:spodRec_pod}
\includegraphics[width=\textwidth]{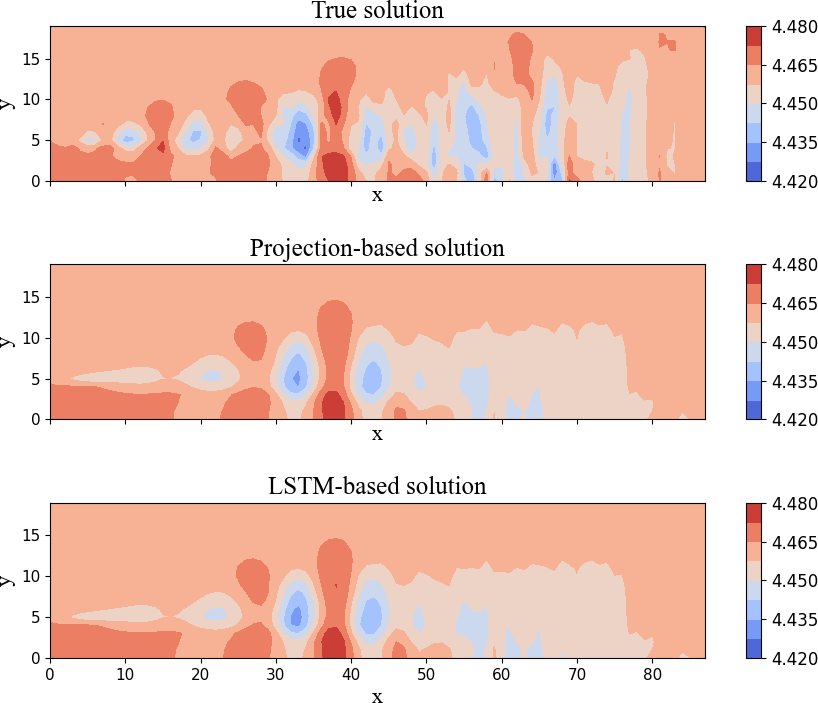}
\caption{POD}
\end{subfigure}
\hspace{0.5cm}%
\begin{subfigure}{.48\textwidth}\label{fig:spodRecs_spod}
\includegraphics[width=\textwidth]{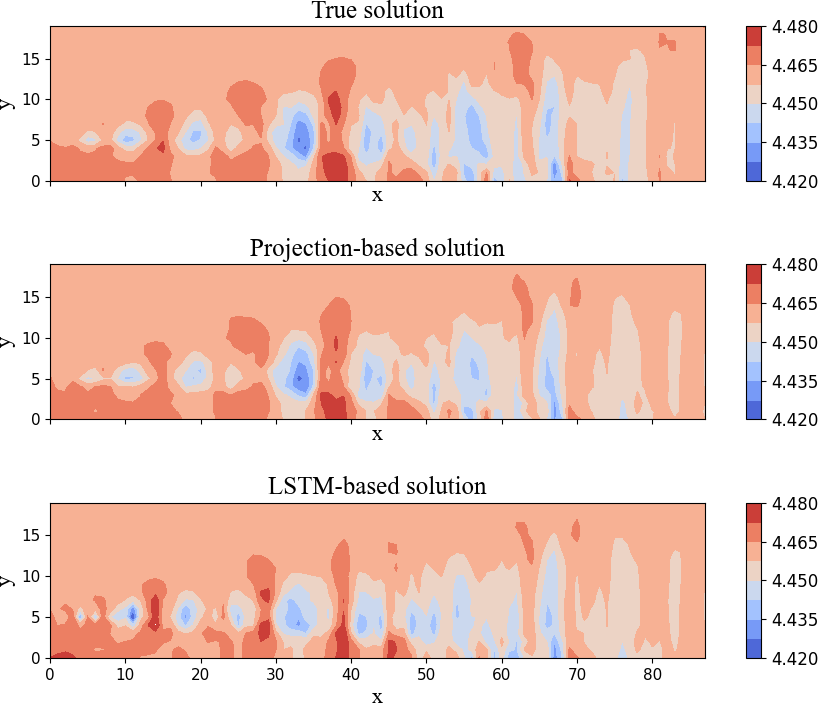}
\caption{SPOD}
\end{subfigure}
\caption{Comparison among: true solution (top), reconstructed solution obtained by projecting the snapshots on the modal basis (middle), reconstructed solution with the predicted LSTM coefficients (bottom). The figures on the left are for POD, while the figures on the right are for SPOD.}
\label{fig:spodRec}
\end{figure}

Figure~\ref{fig:L2vsModes} shows the decay of the averaged L2-norm error against the number of modes for POD (left plot) and SPOD (right plot). We remark again that the SPOD uses all frequencies for the ten modes considered to reconstruct the pressure field. We note how the error of the SPOD reconstruction is generally smaller than the one for the POD reconstruction (note that the y-axis is different for the POD and SPOD). However, the SPOD reconstruction requires a larger amount of data, as we need to consider the frequency dimension for each mode (i.e., 33 frequencies for this test case), that is not present for the POD reconstruction. Therefore, the data compression achieved for SPOD differs from the one achieved for POD, by a factor associated to the number of frequencies retained for the reconstruction, that is: $r_{c} = \text{no.\ modes SPOD} / (\text{no.\ modes POD} \times N_{\text{fr}}) = 1 / N_{\text{fr}})$. 
The neural-network prediction of the coefficients is depicted in figure~\ref{fig:SPODcoeff}, along with the training and validation (that in this case coincides with testing as no hyperparameter optimization was performed) loss as a function of the number of epochs. The results were obtained using the LSTM configuration described in section~\ref{sec:neuralNetworkConfig}, where the single neural network layer employed contained 25 neurons. The total number of epochs for training was equal to 130. The number of past values adopted to make the predictions was $n_{\tau P} = 60$, and the number of future values predicted was $n_{\tau F} = 1$, with a total number of samples for training equal to $n^{(\text{train})}_{\text{samples}} = 749$. The total number of features for SPOD was $n^{(\text{train, SPOD})}_{\text{features}} = 330$ since the first 10 modes were used, resulting in a data compression of $c = 330 / (88 \times 22) = 0.17$, that is, the latent space had a dimension that was only 17\% of the original data.  The total number of features for POD was instead $n^{(\text{train, POD})}_{\text{features}} = 10$, achieving a data compression of $c = 10 / (88 \times 22) = 0.005$, corresponding to a latent space containing an amount of data equal to only  0.5\% of the original data. To this end, the data compression achieved via POD is larger. However, the granularity provided by the SPOD in terms of frequency space is not available. Therefore, the trade-off is application dependent, and revolves around what the frequency information content can provide in terms of analysis and prediction.
Figure~\ref{fig:SPODcoeff} shows how the LSTM prediction for POD (top plots) seems to achieve better performance than the SPOD one (bottom plots). This can be explained by looking at the POD and SPOD coefficients. On the one hand, the POD coefficients are smoother, and resemble a harmonic behaviour. On the other hand, the SPOD coefficients are less smooth, possibly less stationary, and more numerous. A comparison of the first two POD modes (left plot) and the first two SPOD modes (right plot) is depicted in figure~\ref{fig:PODvsSPODcoeff}, where the characteristics described are particularly evident. These different characteristics of the POD and SPOD latent spaces lead to a more challenging neural-network learning task for the SPOD coefficients, thereby negatively affecting the SPOD-LSTM forecasts. One should note that we did not carry out an hyperparameter optimization, nor did we use different neural-network architectures (including attention-mask-based neural network architectures, that can improve the learning performance~\cite{vaswani2017attention}) as these two aspects are out of the scope of this paper. However, these two points can improve the SPOD-LSTM predictions. In addition, having longer time histories for the coefficients may further improve the predictions. 
Figure~\ref{fig:spodRec} shows the comparison among the true solution, the solution reconstructed with the coefficients obtained with a direct projection of the snapshots onto the reduced space, and the one with coefficients predicted with the LSTM network.
The different errors, described in section~\ref{sec:spod_reconstruction}, are reported in table~\ref{tab:pod_spod_errors}. These are the projection error, due to the projection of the true dynamics onto a reduced space, the learning error, due to the neural-network, and the total error, that is a combination of the previous two errors. For the POD, the projection error is consistently larger than the learning error. This indicates that the dominant part of the total error is related to the projection phase, while the LSTM network was very effective in predicting the future values of the coefficients. For the SPOD, the situation is reversed, with the learning error being dominant. Two main elements can explain this different behaviour. The first one is due to the fact that the time series of the SPOD coefficients are more irregular and unsteady with respect to the POD ones and hence more complex to predict (depicted qualitatively in figure~\ref{fig:PODvsSPODcoeff}). The second one is connected to the latent SPOD space dimension, which is significantly larger than the POD one and therefore more prone to accumulation of learning errors. The higher number of features impacts negatively the costs of the training phase of the LSTM network and the one of the model inference. 
This negative impact of the number of modes on the error is confirmed by observing the $L_{\inf}$-norm errors reported in Table \ref{tab:pod_spod_errors}: the learning error is significant when a high number of SPOD modes are taken into account even if the $L_2$-norm error is relatively small, thus indicating that the error is high on a limited number of computational points.

\begin{figure}[H]
\centering
\begin{subfigure}{.49\textwidth}\label{fig:PODvsSPODerr_pod}
\includegraphics[width=\linewidth]{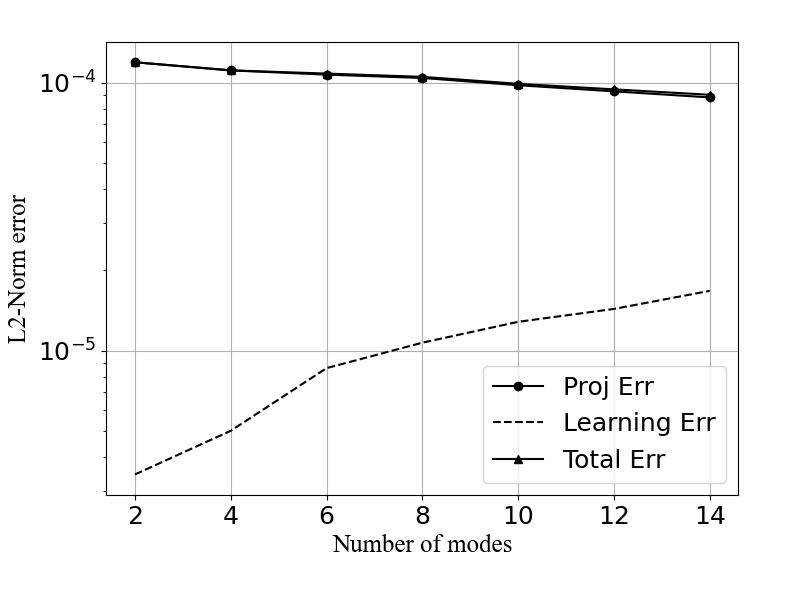}
\caption{POD}
\end{subfigure}
\begin{subfigure}{.49\textwidth}\label{fig:PODvsSPODerr_spod}
\includegraphics[width=\linewidth]{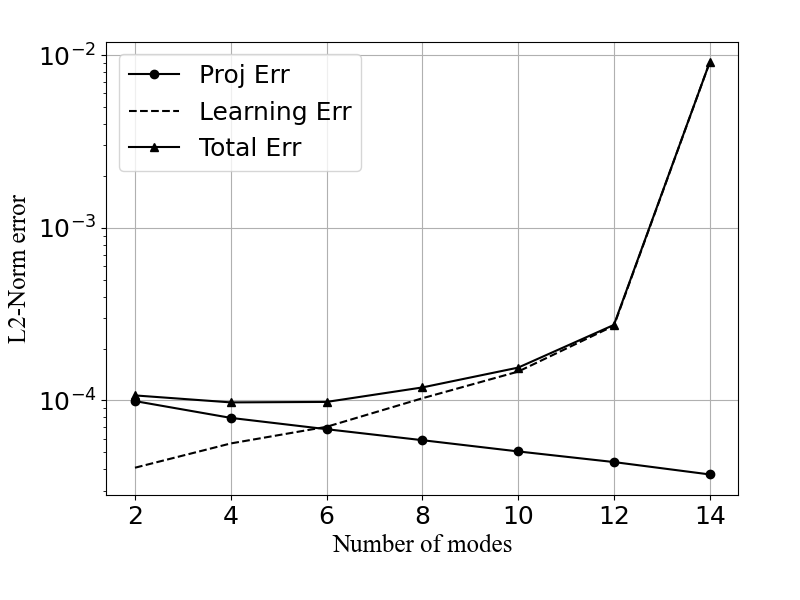}
\caption{SPOD}
\end{subfigure}
\caption{Projection error, learning error, and total error for POD (left) and SPOD (right). Projection and total error for POD are overlapped.}
\label{fig:PODvsSPODerr}
\end{figure}

For both POD and SPOD, the learning error tends to increase with the number of modes, while the projection error decreases (as also shown in figure~\ref{fig:PODvsSPODerr}). 
\begin{table}[H]
\centering
{
\scriptsize
\begin{tabular}{?p{0.8cm}?p{1.1cm}|p{1.1cm}|p{1.1cm} ?p{1.1cm}|p{1.1cm}|p{1.1cm}?p{1.1cm}|p{1.1cm}|p{1.1cm}?}
 \hline
\textbf{POD} &\multicolumn{3}{c?}{\textbf{Projection Error}} &\multicolumn{3}{c?}{\textbf{Learning Error}} &\multicolumn{3}{c?}{\textbf{Total Error}} \\
 \hline
 \hline
 Nr.\ of Modes & $L_{1}$ & $L_2$ & $L_{\infty}$  & $L_1$ & $L_2$ & $L_{\infty}$  & $L_1$ & $L_2$ & $L_{\infty}$ \\
 \hline
 2 & 3.04e-03  &  1.19e-04  &  2.96e-02  &  6.36e-05  &  3.45e-06  &  8.12e-04  &  3.04e-03  &  1.19e-04  &  2.97e-02 \\
 \hline
 4 & 2.81e-03  &  1.11e-04  &  2.88e-02  &  9.95e-05  &  5.02e-06  &  1.13e-03  &  2.81e-03  &  1.11e-04  &  2.88e-02 \\
\hline
 6 & 2.72e-03  &  1.07e-04  &  2.82e-02  &  1.88e-04  &  8.60e-06  &  1.82e-03  &  2.74e-03  &  1.08e-04  &  2.82e-02\\
 \hline
 8 & 2.65e-03  &  1.04e-04  &  2.77e-02  &  2.40e-04  &  1.07e-05  &  2.29e-03  &  2.68e-03  &  1.05e-04  &  2.78e-02 \\
 \hline
 10 & 2.49e-03  &  9.77e-05  &  2.70e-02  &  3.05e-04  &  1.28e-05  &  2.71e-03  &  2.53e-03  &  9.89e-05  &  2.70e-02
 \\
 \hline
 12 & 2.36e-03  &  9.27e-05  &  2.69e-02  &  3.42e-04  &  1.43e-05  &  3.18e-03  &  2.40e-03  &  9.43e-05  &  2.71e-02 \\
 \hline
 14 & 2.26e-03  &  8.80e-05  &  2.64e-02  &  4.06e-04  &  1.67e-05  &  3.69e-03  &  2.31e-03  &  9.00e-05  &  2.67e-02 \\
 \hline
\end{tabular}

\begin{tabular}{?p{0.8cm}?p{1.1cm}|p{1.1cm}|p{1.1cm} ?p{1.1cm}|p{1.1cm}|p{1.1cm}?p{1.1cm}|p{1.1cm}|p{1.1cm}?}
 \hline
 \textbf{SPOD}&\multicolumn{3}{c?}{\textbf{Projection Error}} &\multicolumn{3}{c?}{\textbf{Learning Error}} &\multicolumn{3}{c?}{\textbf{Total Error}} \\
 \hline
 \hline
 Nr.\ of Modes & $L_{1}$ & $L_2$ & $L_{\infty}$  & $L_1$ & $L_2$ & $L_{\infty}$  & $L_1$ & $L_2$ & $L_{\infty}$ \\
 \hline
 2 & 2.81e-03  &  9.92e-05  &  2.32e-02  &  9.50e-04  &  4.08e-05  &  1.21e-02  &  3.01e-03  &  1.07e-04  &  2.52e-02 \\
 \hline
 4 & 2.33e-03  &  7.93e-05  &  1.90e-02  &  1.34e-03  &  5.64e-05  &  1.52e-02  &  2.79e-03  &  9.75e-05  &  2.24e-02 \\
\hline
 6 & 2.04e-03  &  6.82e-05  &  1.45e-02  &  1.69e-03  &  7.04e-05  &  1.86e-02  &  2.81e-03  &  9.81e-05  &  2.17e-02\\
 \hline
 8 &1.79e-03  &  5.89e-05  &  1.18e-02  &  2.45e-03  &  1.03e-04  &  2.55e-02  &  3.26e-03  &  1.19e-04  &  2.65e-02 \\
 \hline
 10 &1.55e-03  &  5.07e-05  &  9.82e-03  &  3.47e-03  &  1.47e-04  &  3.37e-02  &  4.05e-03  &  1.55e-04  &  3.39e-02 \\
 \hline
12 &1.35e-03  &  4.40e-05  &  8.56e-03  &  6.34e-03  &  2.71e-04  &  6.45e-02  &  6.71e-03  &  2.74e-04  &  6.47e-02 \\
\hline
14&1.14e-03  &  3.73e-05  &  7.32e-03  &  1.97e-01  &  9.09e-03  &  2.45e+00  &  1.97e-01  &  9.09e-03  &  2.45e+00 \\
\hline
\end{tabular}
}
\caption{$L_1$, $L_2$, and $L_{\inf}$ norm errors evaluated considering a different number of modes for both POD (top table) and SPOD (bottom table).}
\label{tab:pod_spod_errors}
\end{table}
The number of modes strongly influences the computational time when SPOD is considered since, as described in Section\ref{sec:neuralNetworkConfig}, an increasing number of neural networks has to be trained; Table \ref{tab:computational_time} reports the computational times for both POD and SPOD when a different number of modes is considered. Two different phases are considered in the table aforementioned: the Offline phase, which includes both the generation of the latent space and the training of the neural networks, and the Online phase, during which the high dimensional flow fields are predicted. 
\begin{table}[H]
\centering
{
\scriptsize
\begin{tabular}{?p{2cm}?p{1.1cm}|p{1.1cm}|p{1.1cm} ?p{1.1cm}|}
 \hline
&\multicolumn{2}{c?}{\textbf{POD}} &\multicolumn{2}{c?}{\textbf{SPOD}} \\
 \hline
 Nr. of Modes &   Offline &  Online & Offline  & Online \\
 \hline
 2  & 68.45 & 3.92  &  268.57 & 14.02 \\
 \hline
 4  & 69.29 & 3.82 &   537.73 &  28.63; \\
\hline
 8 & 71.1 &  3.95 & 1166.58 &  56.20; \\
\hline
 10 & 68.08 &  4.04 & 1369.86 & 69.31 \\
 \hline
\end{tabular}
}
\caption{Computational times (in seconds) for both POD and SPOD for a different number of modes.}
\label{tab:computational_time}
\end{table}

\begin{figure}[H]
\centering
\begin{subfigure}{.49\textwidth}\label{fig:POD_timeHorizon}
\includegraphics[width=\linewidth]{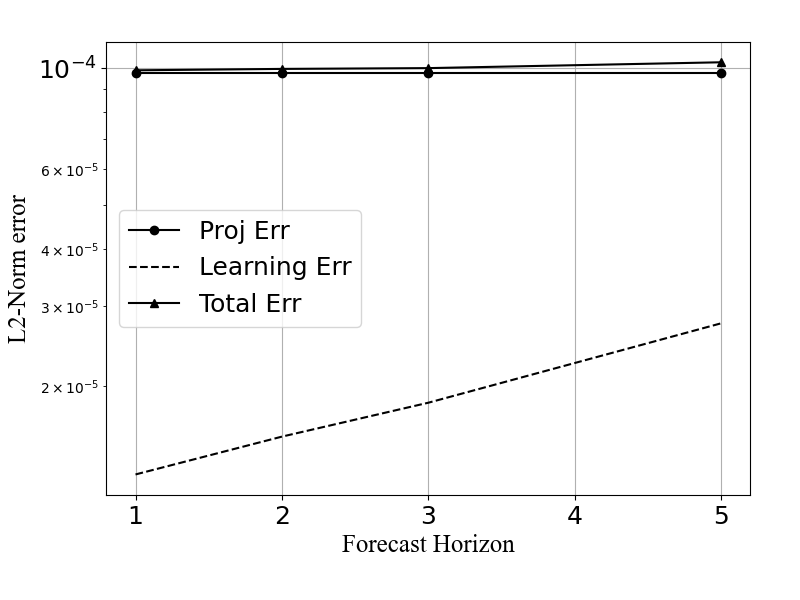}
\caption{POD}
\end{subfigure}
\begin{subfigure}{.49\textwidth}\label{fig:SPOD_timeHorizon}
\includegraphics[width=\linewidth]{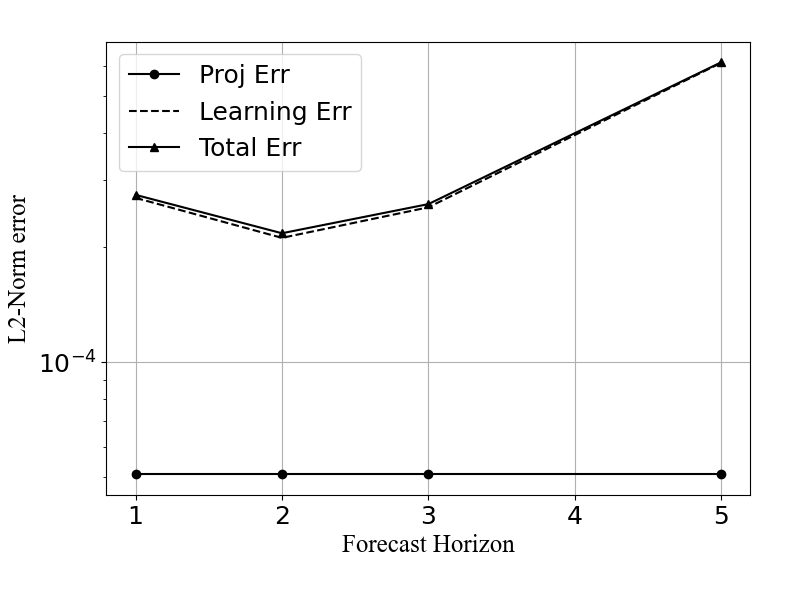}
\caption{SPOD}
\end{subfigure}
\caption{Projection error, learning error, and total error for POD and SPOD using different values of $n_{\tau F}$.}
\label{fig:PODvsSPOD_timeHorizon}
\end{figure}

Figure~\ref{fig:PODvsSPOD_timeHorizon} shows the projection, learning and total error for a five different forecast horizons $n_{\tau F} = [1,2,3,4,5]$. The projection error is obviously constant, while the learning error increases as the forecast horizon increases. This is expected, as a longer forecast horizon is typically more difficult to predict than a shorter one. We can also observe how the learning error drives the total error for SPOD, while the projection error drives the total error for POD. 

\subsubsection{SPOD emulation: sensitivity to frequency}\label{sec:spod_freq_sensitivity}
\begin{figure}[H]
\centering
\begin{subfigure}{.49\textwidth}\label{fig:SPOD_redFreq_1}
\includegraphics[width=\linewidth]{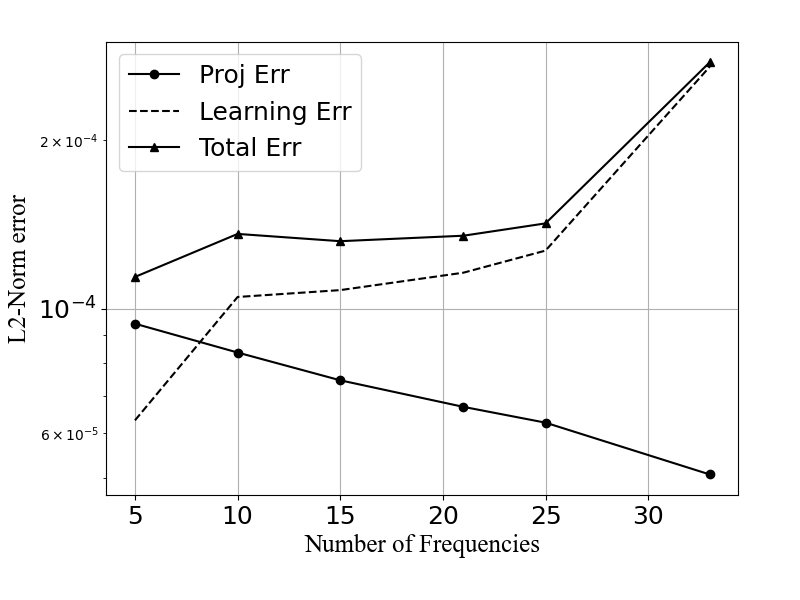}
\caption{$n_{\tau F}=1$}
\end{subfigure}
\begin{subfigure}{.49\textwidth}\label{fig:SPOD_redFreq_3}
\includegraphics[width=\linewidth]{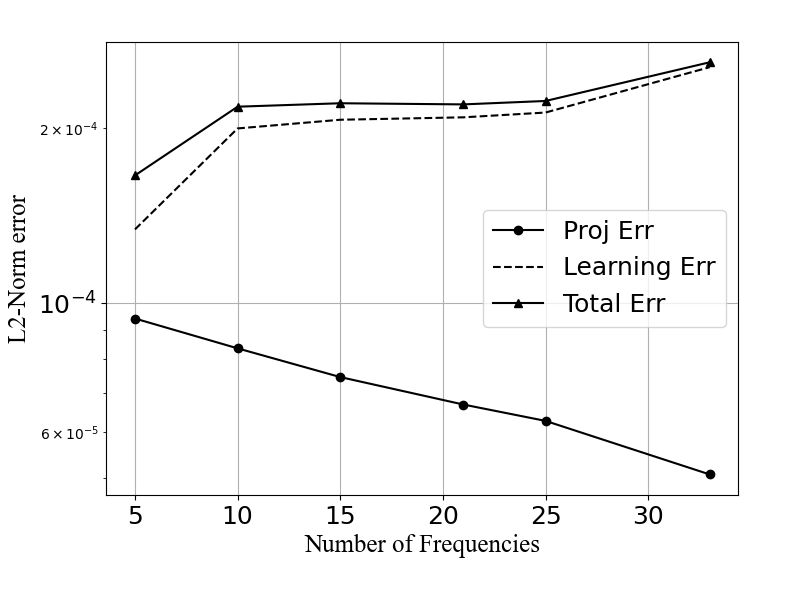}
\caption{$n_{\tau F}=3$}
\end{subfigure}
\caption{Sensitivity to frequency: L2-norm error plot for different numbers of frequency considered and fixed M=10: $n_{\tau F}=1$ (left) and $n_{\tau F} = 3$ (right) }
\label{fig:SPOD_redFreq}
\end{figure}

%
% \begin{figure}[H]
% \centering
% \includegraphics[width=.7\textwidth]{figures/SPOD_redFreq.png}
% \caption{Sensitivity to frequency: L2-norm error plot for different numbers of frequency considered and fixed M=10}
% \label{fig:PODreconstructed}
% \end{figure}
In this section, we show the sensitivity of the emulation framework to the number of frequencies retained for each SPOD mode. In particular, figure~\ref{fig:SPOD_redFreq} shows the projection, learning and total errors for forecast horizon $n_{\tau F} = 1$ (left plot), and $n_{\tau F} = 3$ (right plot). As we increase the number of frequencies retained, the learning error increases, while the projection error decreases, as expected. Table~\ref{tab:spod_errors_freqs} tabulated the errors depicted in figure~\ref{fig:SPOD_redFreq}. Here, it is important to note that by removing higher frequencies we have a smaller learning error. This aspect can have important applications when we need to perform predictive tasks on noisy data. 

\begin{table}[H]
\centering
{
\scriptsize
\begin{tabular}{?p{0.8cm}?p{1.1cm}|p{1.1cm}|p{1.1cm} ?p{1.1cm}|p{1.1cm}|p{1.1cm}?p{1.1cm}|p{1.1cm}|p{1.1cm}?}
 \hline
\textbf{SPOD} &\multicolumn{3}{c?}{\textbf{Projection Error}} &\multicolumn{3}{c?}{\textbf{Learning Error}} &\multicolumn{3}{c?}{\textbf{Total Error}} \\
 \hline
 \hline
 Nr. of Freqs &  $L_{1}$ & $L_2$ & $L_{\infty}$  & $L_1$ & $L_2$ & $L_{\infty}$  & $L_1$ & $L_2$ & $L_{\infty}$ \\
 \hline
 5 & 2.70e-03  &  9.41e-05  &  2.11e-02  &  1.50e-03  &  6.33e-05  &  1.11e-02  &  3.19e-03  &  1.14e-04  &  2.36e-02 \\
 \hline
 10 & 2.44e-03 &  8.36e-05  & 1.74e-02  &  2.48e-03 &   1.05e-04  &  2.07e-02 &   3.70e-03  &  1.36e-04  &  2.75e-02 \\
\hline
 15 & 2.21e-03  &  7.46e-05  &  1.48e-02  &  2.55e-03  &  1.08e-04  &  2.30e-02  &  3.59e-03  &  1.32e-04  &  2.76e-02\\
 \hline
 21 & 2.02e-03  &  6.69e-05  &  1.27e-02  &  2.73e-03  &  1.16e-04  &  2.59e-02  &  3.65e-03  &  1.35e-04  &  2.85e-02 \\
 \hline
 25 & 1.91e-03  &  6.27e-05  &  1.19e-02  &  3.01e-03  &  1.27e-04  &  2.98e-02  &  3.82e-03  &  1.42e-04  &  3.15e-02 \\
 \hline
 33 & 1.55e-03  &  5.07e-05  &  9.82e-03  &  6.37e-03  &  2.70e-04  &  5.77e-02  &  6.82e-03  &  2.75e-04  &  5.81e-02 \\
 \hline
\end{tabular}
}
\caption{$L_1$, $L_2$, and $L_{inf}$ norm errors evaluated retaining a different number of SPOD frequencies.}
\label{tab:spod_errors_freqs}
\end{table}

\subsection{Frequency domain reconstruction}\label{sec:spod_freq_reconstruction_sensitivity}

In Section~\ref{sec:spod_latent} two distinct techniques for the evaluation of the expansion coefficients were presented; the aim of this section is to provide a comparison between the frequency-domain approach and the time-domain one used in the previous sections.

\begin{figure}[H]
\centering
\begin{subfigure}{.49\textwidth}
\includegraphics[width=\linewidth]{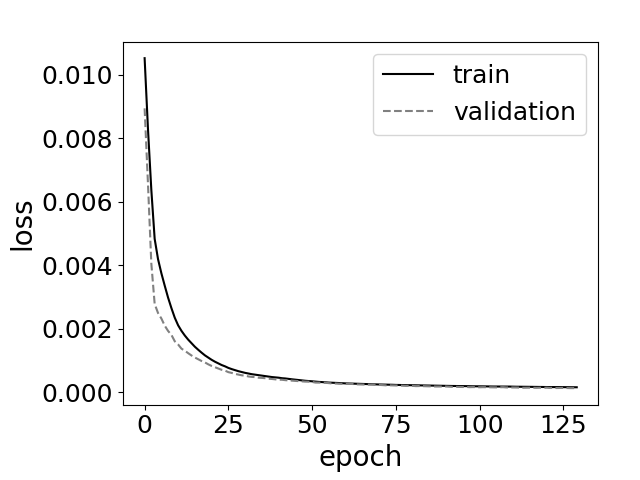}
\caption{Training history}
\label{fig:jet_freqDomain_history}
\end{subfigure}
\begin{subfigure}{.49\textwidth}
\includegraphics[width=\linewidth]{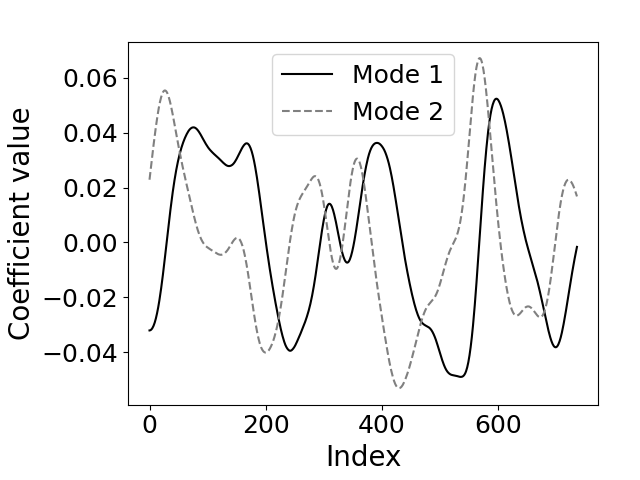}
\caption{Modes}
\label{fig:jet_freqDomain_coeffs}
\end{subfigure}
\caption{Loss histories of the training phase (on the left) and comparison between the first two modes series used for the training of the LSTM network.}
\label{fig:SPODcoeff_freqDomain}
\end{figure}

One of the main drawback of the frequency-domain approach, when applied to neural network, is that the dimension of the latent space is proportional to the number of blocks $L$ and not to the number of time snapshots. Therefore, in order to have an adequate number of samples for training the LSTM-based neural network, the data set containing the 800 flow realizations of the training set has been split into 687 blocks with an overlap of $98\%$, i.e. two contiguous blocks differ for just one snapshot. 
The training phase of the $L\times N_f$ expansion coefficients has been carried out with the aid of the neural network configured as described in Section \ref{sec:neuralNetworkConfig}. The histories of the losses recorded during the training are reported in Figure \ref{fig:jet_freqDomain_history} and they demonstrate the suitability of the chosen hyper-parameters also for the current task.

The coefficients associated to the first and the second modes are shown in Figure \ref{fig:jet_freqDomain_coeffs} and their shape is much more regular and sinusoidal-like with respect to the ones obtained in the time-domain (Figure \ref{fig:PODvsSPODcoeff}). 
Following this idea, we expect a lower learning error for the frequency domain case and this is confirmed by the results reported in Table \ref{tab:spod_errors_modes_freq_domain} where norm errors are reported when different numbers of modes are considered.
On the other hand the oblique projection used in the time-domain approach guarantees the best approximation in a least-square sense and this provide a justification for its smaller projection errors that can be observed.
In the present work we have chosen to work mainly with the time domain reconstruction since it leads to a number of training samples which is the same of the POD case and this guarantees a consistent comparison between the two techniques.

\begin{table}[H]
\centering
{
\scriptsize
\begin{tabular}{?p{0.8cm}?p{1.1cm}|p{1.1cm}|p{1.1cm} ?p{1.1cm}|p{1.1cm}|p{1.1cm}?p{1.1cm}|p{1.1cm}|p{1.1cm}?}
 \hline
\textbf{SPODf} &\multicolumn{3}{c?}{\textbf{Projection Error}} &\multicolumn{3}{c?}{\textbf{Learning Error}} &\multicolumn{3}{c?}{\textbf{Total Error}} \\
 \hline
 \hline
 Nr. of Freqs &  $L_{1}$ & $L_2$ & $L_{\infty}$  & $L_1$ & $L_2$ & $L_{\infty}$  & $L_1$ & $L_2$ & $L_{\infty}$ \\
 \hline
2 & 3.65e-03 & 1.34e-04 & 2.79e-02 & 2.34e-04 & 1.11e-05 & 4.21e-03 & 3.65e-03 & 1.34e-04 & 2.79e-02 \\
 \hline
4 & 3.58e-03 & 1.28e-04 & 2.76e-02 & 2.60e-04 & 1.21e-05 & 4.58e-03 & 3.58e-03 & 1.28e-04 & 2.74e-02 \\
 \hline
6 & 3.53e-03 & 1.27e-04 & 2.84e-02 & 2.77e-04 & 1.27e-05 & 4.72e-03 & 3.53e-03 & 1.27e-04 & 2.83e-02 \\
 \hline
8 & 3.54e-03 & 1.25e-04 & 2.64e-02 & 2.91e-04 & 1.31e-05 & 4.67e-03 & 3.54e-03 & 1.25e-04 & 2.63e-02 \\
 \hline
10 & 3.53e-03 & 1.25e-04 & 2.72e-02 & 3.01e-04 & 1.33e-05 & 4.67e-03 &  3.53e-03 & 1.25e-04 & 2.70e-02 \\
 \hline
12 & 3.52e-03 & 1.25e-04 & 2.71e-02 & 3.12e-04 & 1.36e-05 & 4.67e-03 & 3.53e-03 & 1.25e-04 & 2.71e-02 \\
 \hline
14 & 3.50e-03 & 1.24e-04 &  2.65e-02 & 3.19e-04 & 1.38e-05 & 4.70e-03 & 3.51e-03 & 1.25e-04 & 2.64e-02 \\
 \hline
\end{tabular}
}
\caption{$L_1$, $L_2$, and $L_{inf}$ norm errors evaluated retaining a different number of SPOD modes for frequency domain reconstruction.}
\label{tab:spod_errors_modes_freq_domain}
\end{table}

\subsection{Geophysical problems: the Madden-Julian Oscillation}\label{sec:results_mjo}
The second test case which was performed considering a geophysical tropical phenomenon known as the Madden-Juklian Oscillation (MJO). This is an intraseasonal phenomenon whose characteristic period varies between $30$ and $90$ days and it is the result of a large-scale coupling between atmospheric circulation and deep convection. This pattern slowly propagates eastward with a speed of $4$ to $8 \,m\,{s^{-1}}$. Although it is commonly known as an oscillation, the large ranges of variation of both its local period (30-90 days) and its speed of propagation suggest that MJO is a rather irregular phenomenon which can be seen at a large-scale level as a mix of multiple high-frequency, small-scale convective phenomena. Further information about MJO can be found in~\cite{zhang2005madden}.

% The Madden-Julian Oscillation (MJO) is an intraseasonal phenomenon that characterizes the tropical atmosphere. It has a characteristic period which varies between $30$ and $90$ days and it is the result of a large-scale coupling between atmospheric circulation and deep convection. This pattern slowly propagates eastward with a speed of $4$ to $8 \,m\,{s^{-1}}$. Although it is commonly known as an oscillation, the large ranges of variation of both its local period (30-90 days) and its speed of propagation suggest that MJO is a rather irregular phenomenon. This implies that the MJO can be seen at a large-scale level as a mix of multiple high-frequency, small-scale convective phenomena. Further information about MJO can be found in~\cite{zhang2005madden}.

The aim of the present case is to test the capabilities of the SPOD-based approach in isolating and predicting the evolution in time of phenomena characterized by specific subsets of frequencies. The capability of breaking down complex flow fields in simpler structures united by similar characteristic  time-scales could be beneficial for applications involving, for instance, data denoising, filtering, and frequency analysis. 

In this test case the starting dataset was composed by 5000 time snapshots, each of these providing a representation of the total precipitation with a resolution of 12 hours. The training data set consisted of the first $75\%$ samples (equal to  $n^{(\text{train})}_{\text{samples}}=3750$), while the remaining $25\%$ samples (equal to $n^{(\text{test})}_{\text{samples}}=1250$) were used for testing and validation purposes. The spatial dimension of the data consists of $480 \times 241 = 115680$ points, where $480$ are longitude points and $241$ are latitude points. 

% The present test case is proposed with the aim of showing the capabilities of the SPOD-based approach in isolating and predicting the evolution of phenomena characterized by a specific range of frequencies; this can have interesting applications for different purposes, for instance for data denoising, filtering, and frequency analysis. In this test case the starting dataset was composed by 5000 snapshots representing the total precipitation with a resolution of 12 hours. The training data consisted of the first $75\%$ samples (equal to  $n^{(\text{train})}_{\text{samples}}=3750$), while the remaining $25\%$ samples (equal to $n^{(\text{test})}_{\text{samples}}=1250$) were used for testing and validation purposes. The spatial dimension of the data consists of $480 \times 241 = 115680$ points, where $480$ are longitude points and $241$ are latitude points. 

The snapshots in the training set were split into $L = 5$ statistically independent blocks, each of them containing 730 snapshots with a $\ell_{\text{ovlp}} = 0\%$ overlap. The weighting for the inner product is based on the surface element of a sphere $\text{d} S_r=r^2 \sin\theta \, \text{d}\theta \, \text{d}\varphi$ with radius $r$ over $(\theta,\varphi)\in [-\pi/2,\pi/2]\times[-\pi,\pi]$, where we denote the spherical coordinate triplet $\vb{x}=(r,\theta,\varphi)$ consisting of radius, geographic latitude and geographic longitude, respectively. The choice of parameters just introduced provides up to 5 modes and 366 frequencies. The SPOD modes were generated using the methods illustrated in section \ref{sec:methods}, starting from the snapshots containing the rain precipitation fluctuations computed by subtracting the time mean from each realization. 

The neural network used for learning the dynamics of the SPOD latent space is the same whose configuration has been described in section~\ref{sec:neuralNetworkConfig}.

% To compute the SPOD modes, the snapshots of the training set were subdivided into $L = 5$ statistically independent blocks, each of them containing 730 snapshots with a $\ell_{\text{ovlp}} = 0\%$ overlap. The weighting for the inner product is based on the surface element of a sphere $\text{d} S_r=r^2 \sin\theta \, \text{d}\theta \, \text{d}\varphi$ with radius $r$ over $(\theta,\varphi)\in [-\pi/2,\pi/2]\times[-\pi,\pi]$, where we denote the spherical coordinate triplet $\vb{x}=(r,\theta,\varphi)$ consisting of radius, geographic latitude and geographic longitude, respectively. This choice of parameters provides up to 5 modes and 366 frequencies. The SPOD modes were generated using the methods illustrated in section \ref{sec:methods}, starting from the snapshots containing the rain precipitation fluctuations computed by subtracting the time mean from each realization. 

In the following subsections, we provide a comparison between the results obtained by using the POD-based approach and SPOD-based one, and the sensitivity with respect to some significant parameters such as the number of modes and the number of frequencies retained for the reconstruction is discussed in section~\ref{subsec:mjo_spod_freqs} and section~\ref{subsec:mjo_pod_spod}, respectively.

% In the following two subsections, we show the comparison between POD and SPOD emulation, and the sensitivity with respect to the number of modes (section~\ref{subsec:mjo_spod_freqs}), as well as the sensitivity of SPOD emulation to the frequencies retained (section~\ref{subsec:mjo_pod_spod}).

\subsubsection{POD vs.\ SPOD emulation: Sensitivity to number of modes}\label{subsec:mjo_pod_spod}
In this section the influence on the errors of the numbers of modes is investigated.

\begin{figure}[H]
    \centering
    \begin{subfigure}{.48\textwidth}
        \includegraphics[width=\textwidth]{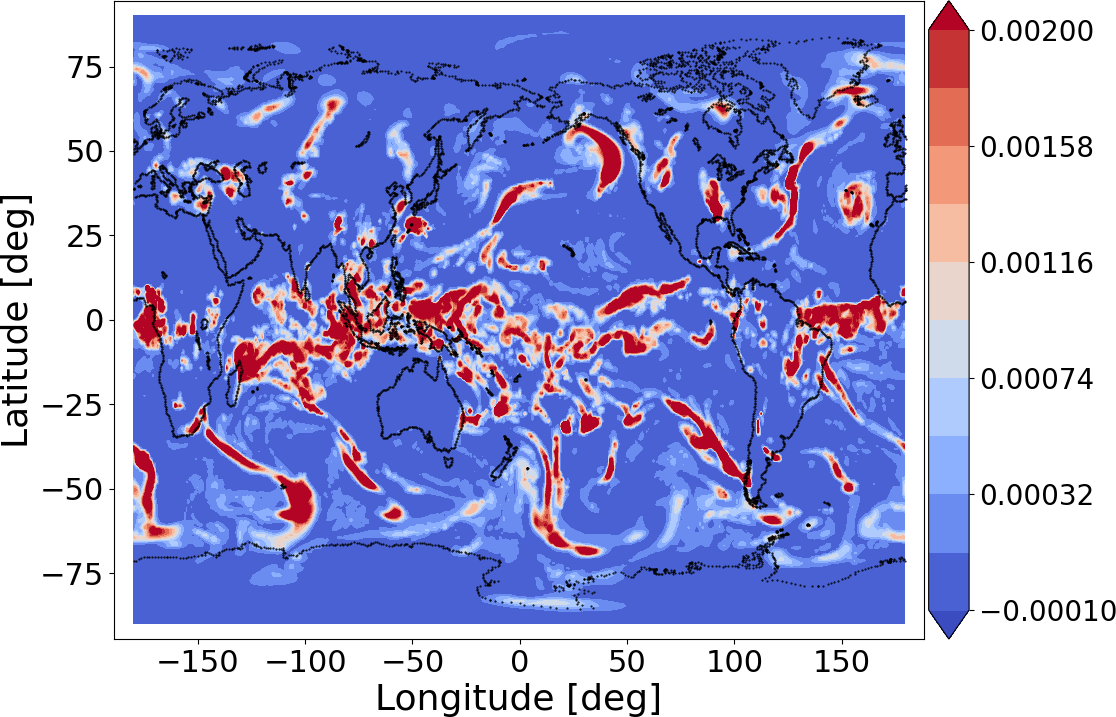}
        \caption{}
        \label{fig:MJO_flowRealization}
    \end{subfigure}
    \hspace{0.5cm}%
    \begin{subfigure}{.48\textwidth}
        \includegraphics[width=\textwidth]{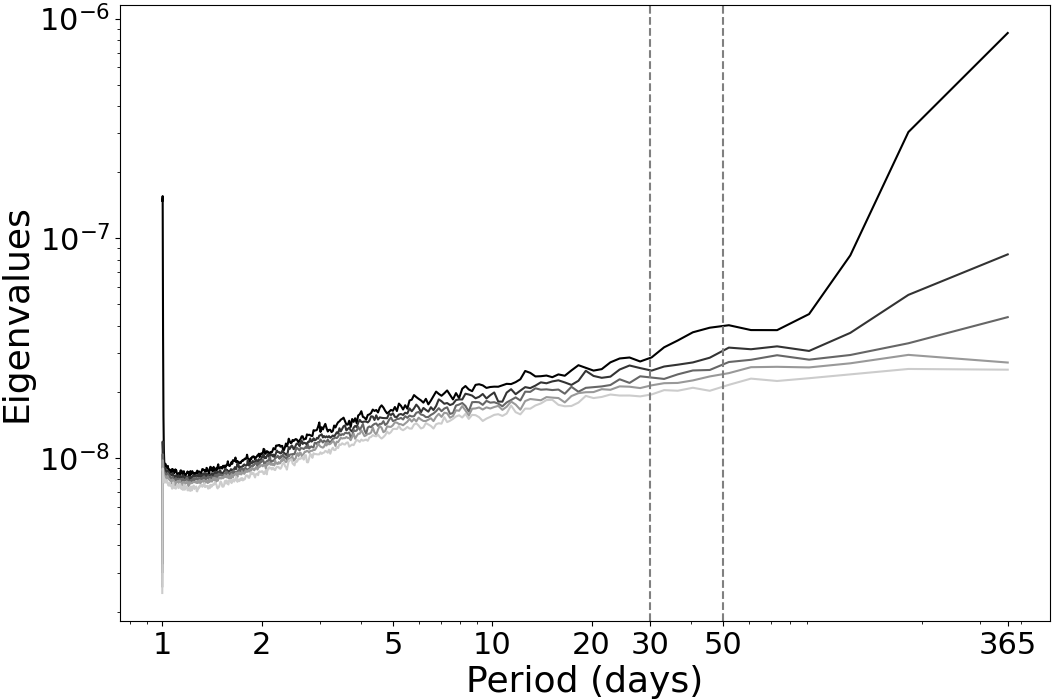}        \caption{}
        \label{fig:MJO_eigsVSperiod8}
    \end{subfigure}
    \caption{The 100-th snapshot of the database: total precipitation of rain (left plot). Eigenvalues of the first 5 modes plotted against period obtained using SPOD (right plot).}
    \label{fig:MJO_fig}
\end{figure}

First, the SPOD modes are evaluated starting from the flow realizations contained in the training dataset; a visualization of a snapshot belonging to this database is reported in figure~\ref{fig:MJO_flowRealization}. The eigenvalues associated to the first five modes at different characteristic periods are represented in the plot shown in figure~\ref{fig:MJO_eigsVSperiod8}. The influence of the first mode is noticeable in the region associated to low frequencies; in the plot dashed lines help to identify the characteristic period in which the MJO phenomenon is expected to be significant. 

% Figure~\ref{fig:MJO_snapshot} shows the 100-th snapshot of the total rain precipitation field (left plot) and the first 5 eigenvalues (associated to the SPOD modes) plotted against period (right plot), where the two vertical dashed lines indicate the approximate period of MJO (30 to 50 days). 
%
\begin{figure}[H]
\centering
\begin{subfigure}{.49\textwidth}\label{fig:MJO_SPOD_all_freqs_1_modes_ntF1}
\includegraphics[width=\linewidth]{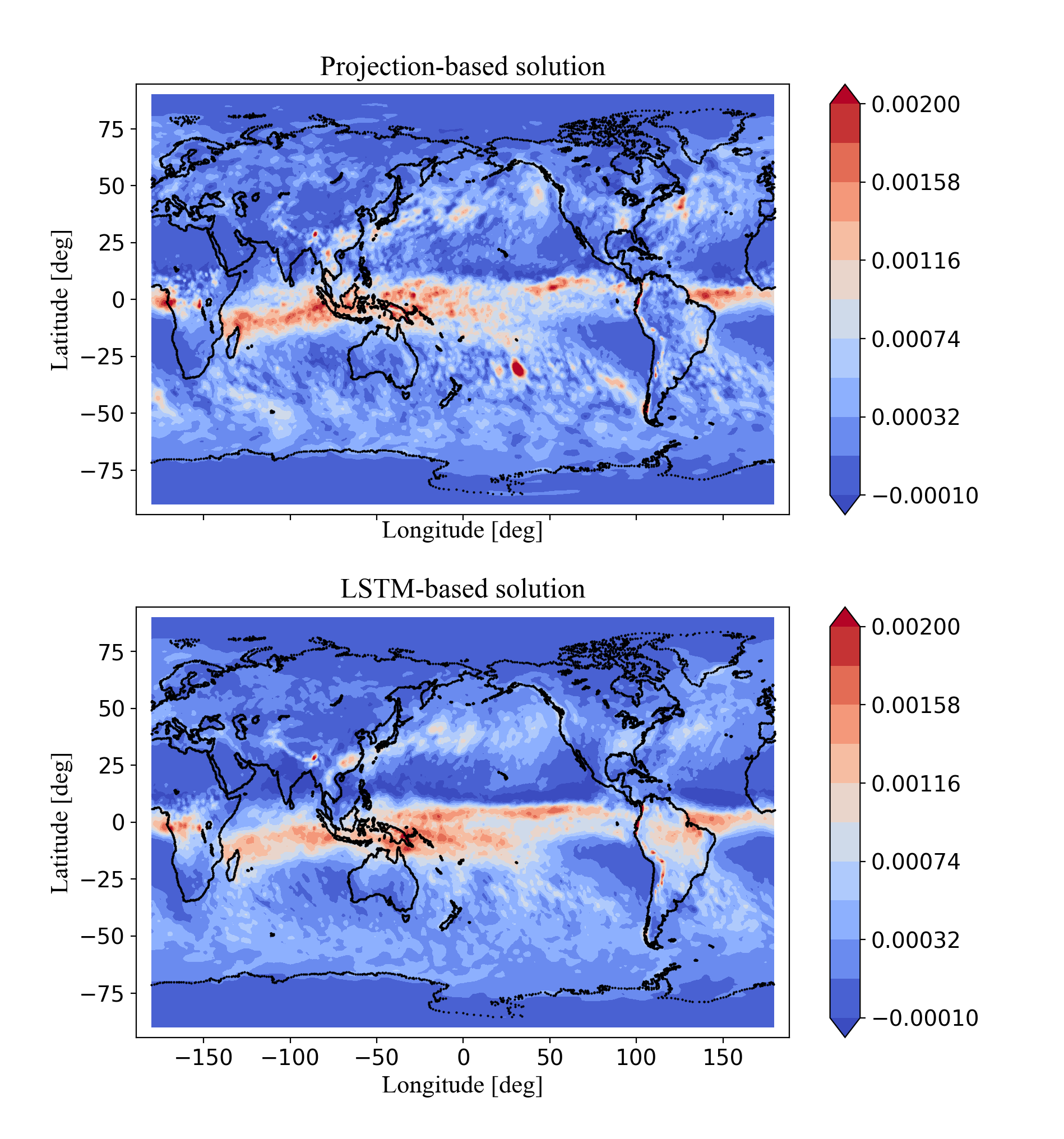}
\caption{1 mode, all frequencies}
\end{subfigure}
\begin{subfigure}{.49\textwidth}\label{fig:MJO_SPOD_all_freqs_5_modes_ntF1}
\includegraphics[width=1.01\linewidth]{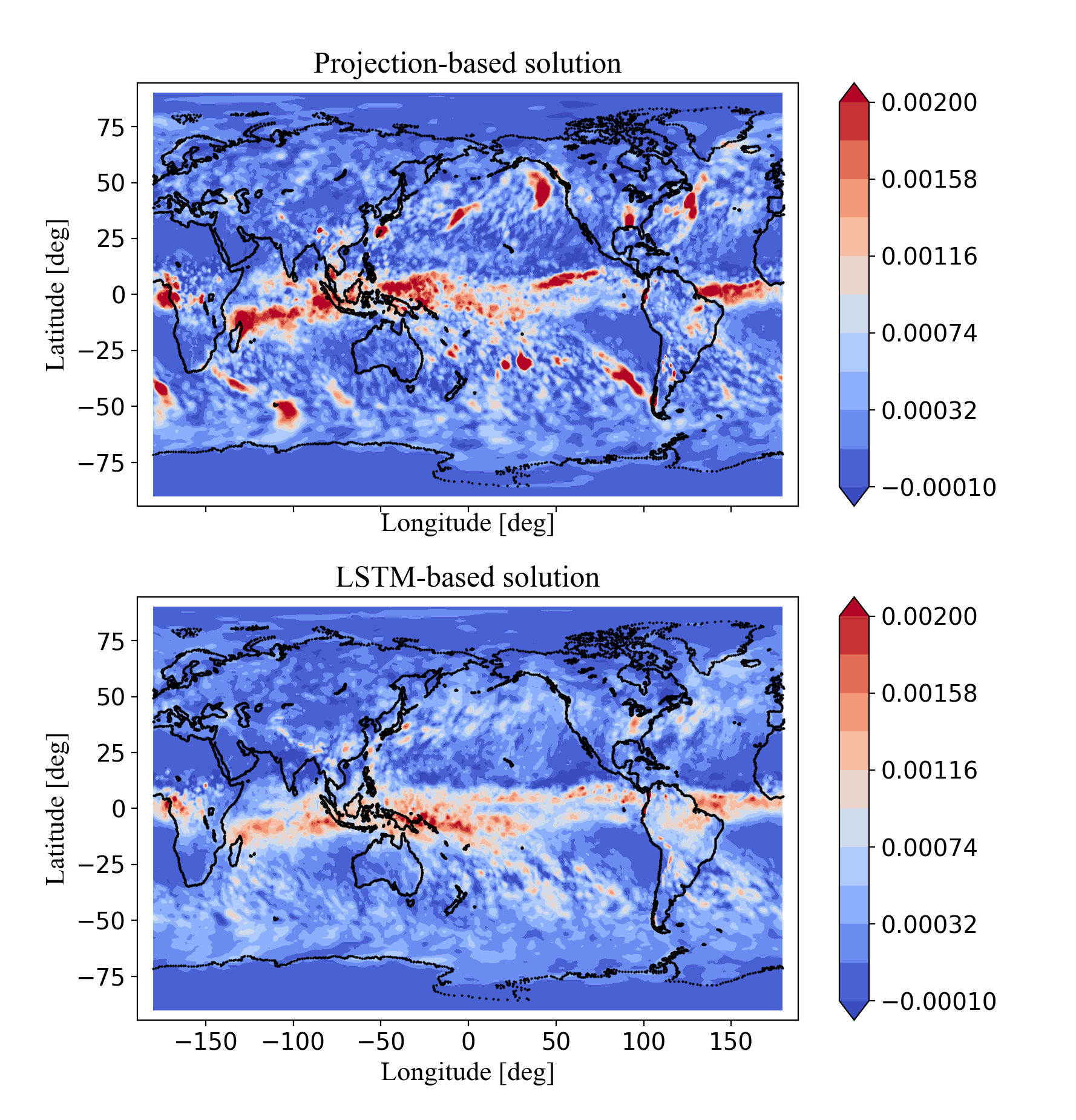}
\caption{5 modes, all frequencies}
\end{subfigure}
\caption{Projection solution (top) vs LSTM prediction (bottom) at the same time snapshot = 100-th, retaining 1 SPOD mode (left) and 5 SPOD modes (right). For each mode retained, all frequencies were used. The forecast horizon adopted was $n_{\tau F} = 1$ (corresponding to 12 hours).}
\label{fig:MJO_SPOD_all_freqs_ntF1}
\end{figure}

Following the procedure described in \ref{sec:methods}, the coefficients are evaluated and then the latent dynamics is learnt via a neural network for both the SPOD-based approach and the POD-based one. The SPOD coefficients are obtained via time-domain reconstruction, as described in section~\ref{sec:time-domain-reconstruction}.
The results of the procedure are shown in figure~\ref{fig:MJO_SPOD_all_freqs_ntF1}, where the reconstructions of the solution in time after projection (top plots) and after LSTM learning (bottom plots) are shown. The left plots were obtained using 1 SPOD mode, while the right plots using all the 5 SPOD modes. In all the plots in figure~\ref{fig:MJO_SPOD_all_freqs_ntF1} all frequencies were retained and a forecast horizon $n_{\tau F} = 1$ was used. By comparing these figures against the snapshot depicted in figure~\ref{fig:MJO_flowRealization} one can observe that, in this specific test case, both the projection- and LSTM-based reconstructions are visibly different from the true solution. Yet, the LSTM-based solution resemble the projection based one. Therefore, the largest source of error in this case is due to the projection, not to the learning task.

%====================
\begin{table}[H]
\centering
{
\scriptsize
\begin{tabular}{?p{0.8cm}?p{1.1cm}|p{1.1cm}|p{1.1cm} ?p{1.1cm}|p{1.1cm}|p{1.1cm}?p{1.1cm}|p{1.1cm}|p{1.1cm}?}
 \hline
\textbf{POD} &\multicolumn{3}{c?}{\textbf{Projection Error}} &\multicolumn{3}{c?}{\textbf{Learning Error}} &\multicolumn{3}{c?}{\textbf{Total Error}} \\
 \hline
 \hline
 Nr.\ of Modes &  $L_{1}$ & $L_2$ & $L_{\infty}$  & $L_1$ & $L_2$ & $L_{\infty}$  & $L_1$ & $L_2$ & $L_{\infty}$ \\
\hline
1  &  4.86e-04  &  2.36e-06  &  3.12e-02  &  7.49e-06  &  3.82e-08  &  9.91e-05  &  4.86e-04  &  2.36e-06  &  3.12e-02 \\
\hline
2  &  4.83e-04  &  2.34e-06  &  3.12e-02  &  1.05e-05  &  6.07e-08  &  2.42e-04  &  4.83e-04  &  2.35e-06  &  3.12e-02 \\
\hline
3 &  4.80e-04  &  2.34e-06  &  3.12e-02  &  1.30e-05  &  7.26e-08  &  3.90e-04  &  4.81e-04  &  2.34e-06  &  3.12e-02 \\
\hline
4  &  4.79e-04  &  2.34e-06  &  3.12e-02  &  1.55e-05  &  8.64e-08  &  4.13e-04  &  4.80e-04  &  2.34e-06  &  3.12e-02 \\
\hline
5  &  4.79e-04  &  2.33e-06  &  3.11e-02  &  3.67e-05  &  2.37e-07  &  1.05e-03  &  4.83e-04  &  2.35e-06  &  3.12e-02 \\
 \hline
\end{tabular}

\begin{tabular}{?p{0.8cm}?p{1.1cm}|p{1.1cm}|p{1.1cm} ?p{1.1cm}|p{1.1cm}|p{1.1cm}?p{1.1cm}|p{1.1cm}|p{1.1cm}?}
 \hline
\textbf{SPOD} &\multicolumn{3}{c?}{\textbf{Projection Error}} &\multicolumn{3}{c?}{\textbf{Learning Error}} &\multicolumn{3}{c?}{\textbf{Total Error}} \\
 \hline
 \hline
 Nr.\ of Modes &  $L_{1}$ & $L_2$ & $L_{\infty}$  & $L_1$ & $L_2$ & $L_{\infty}$  & $L_1$ & $L_2$ & $L_{\infty}$ \\
 \hline
1  &  3.33e-04  &  1.91e-06  &  2.97e-02  &  1.81e-04  &  7.54e-07  &  3.51e-03  &  3.27e-04  &  2.05e-06  &  3.18e-02 \\
\hline
2  &  3.38e-04  &  1.83e-06  &  2.81e-02  &  2.42e-04  &  1.01e-06  &  5.58e-03  &  3.44e-04  &  2.08e-06  &  3.17e-02 \\
\hline
3  &  3.38e-04  &  1.75e-06  &  2.67e-02  &  2.86e-04  &  1.20e-06  &  7.23e-03  &  3.63e-04  &  2.11e-06  &  3.17e-02 \\
\hline
4  &  3.35e-04  &  1.69e-06  &  2.56e-02  &  3.04e-04  &  1.29e-06  &  8.65e-03  &  3.61e-04  &  2.11e-06  &  3.17e-02 \\
\hline
5  &  3.31e-04  &  1.63e-06  &  2.45e-02  &  3.34e-04  &  1.42e-06  &  9.84e-03  &  3.80e-04  &  2.15e-06  &  3.17e-02 \\
 \hline
\end{tabular}
\caption{$L_1$, $L_2$, and $L_{inf}$ norm errors evaluated considering a different number of modes for both POD (top table) and SPOD (bottom table).}
\label{tab:MJO_PODvsSPOD_modes}
}
\end{table}

%====================

\begin{figure}[H]
\centering
\begin{subfigure}{.49\textwidth}\label{fig:MJO_errors_POD_ntF1}
\includegraphics[width=\linewidth]{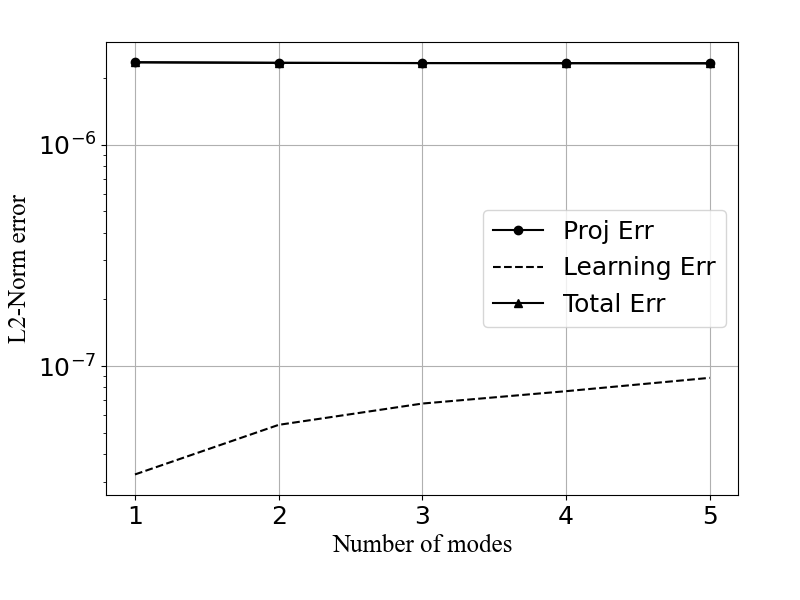}
\caption{POD}
\end{subfigure}
\begin{subfigure}{.49\textwidth}\label{fig:MJO_errors_all_freqs_ntF1}
\includegraphics[width=\linewidth]{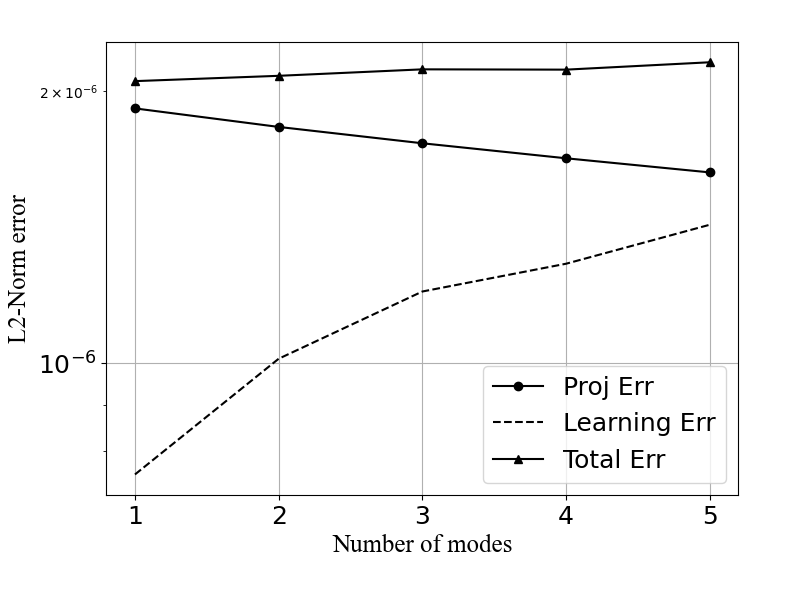}
\caption{SPOD}
\end{subfigure}
\caption{Projection error, learning error, and total error plotted against the number of modes for forecast horizon $n_{\tau F} = 1$ (corresponding to 12 hours): POD (left) and SPOD retaining all frequencies for each mode (right). For POD, projection and total error are overlapped.}
\label{fig:MJO_POD_vs_SPOD_L2_ntF1}
\end{figure}

\begin{figure}[H]
\centering
\begin{subfigure}{.49\textwidth}\label{fig:MJO_errors_POD_ntF3}
\includegraphics[width=\linewidth]{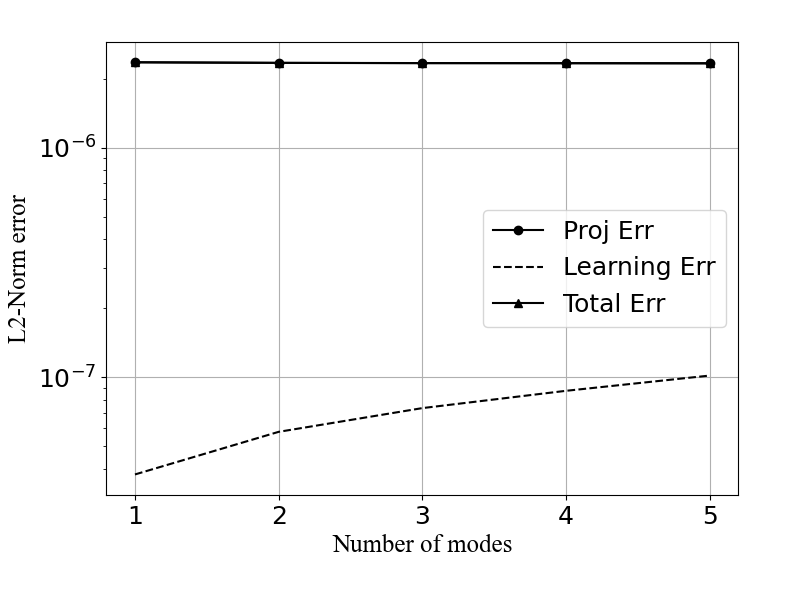}
\caption{POD}
\end{subfigure}
\begin{subfigure}{.49\textwidth}\label{fig:MJO_errors_all_freqs_ntF3}
\includegraphics[width=\linewidth]{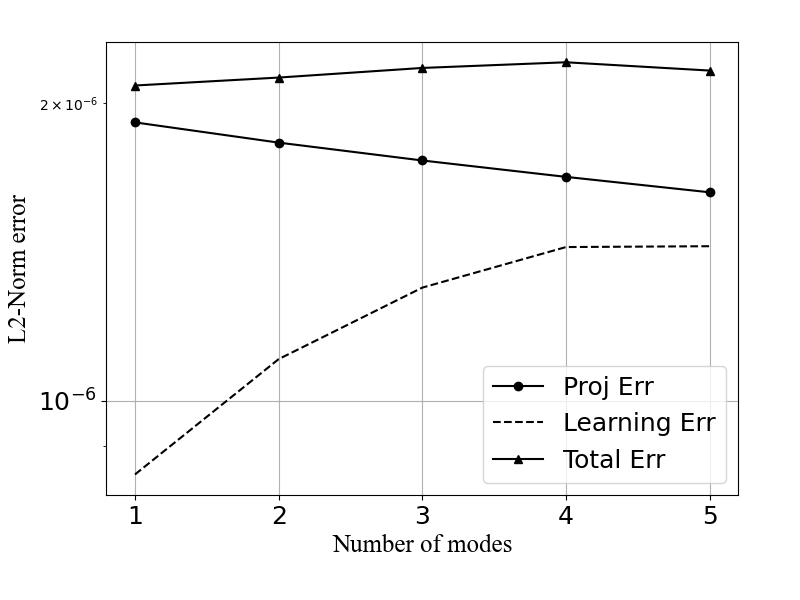}
\caption{SPOD}
\end{subfigure}
\caption{Projection error, learning error, and total error plotted against the number of modes for forecast horizon $n_{\tau F} = 3$ (corresponding to 36 hours): POD (left) and SPOD retaining all frequencies for each mode (right). For POD, projection and total error are overlapped.}
\label{fig:MJO_POD_vs_SPOD_L2_ntF3}
\end{figure}

%
% %
% \begin{figure}[H]
% \centering
% \begin{subfigure}{.49\textwidth}\label{fig:MJO_SPOD_all_freqs_1_modes_ntF3}
% \includegraphics[width=\linewidth]{figures/MJO_SPOD_all_freqs_1_modes_ntF3.png}
% \caption{1 mode, all frequencies}
% \end{subfigure}
% \begin{subfigure}{.49\textwidth}\label{fig:MJO_SPOD_all_freqs_1_modes_ntF3}
% \includegraphics[width=\linewidth]{figures/MJO_SPOD_all_freqs_5_modes_ntF3.png}
% \caption{5 modes, all frequencies}
% \end{subfigure}
% \caption{Projection solution (top) vs LSTM prediction (bottom) at the same time snapshot = 100-th, retaining 1 SPOD mode (left) and 5 SPOD modes (right). For each mode retained, all frequencies were used. The forecast horizon adopted was $n_{\tau F} = 3$ days.}
% \label{fig:MJO_SPOD_all_freqs_ntF3}
% \end{figure}
% %
Figures~\ref{fig:MJO_POD_vs_SPOD_L2_ntF1} and \ref{fig:MJO_POD_vs_SPOD_L2_ntF3} show the projection, learning and total errors as a function of the number of POD (left plots), and SPOD (right plots) modes, for forecast horizons $n_{\tau F} = 1$ (12 hours) and $n_{\tau F} = 3$ (36 hours), respectively. Similarly to the test case presented in section~\ref{sec:comprJet}, increasing the number of modes for both POD and SPOD increases the learning error, while there is a slight decrease in the projection error for both POD and SPOD. The total error for POD is driven by the projection error, while the total error for SPOD is a combination of the projection and learning errors. The plots for SPOD in both figure~\ref{fig:MJO_POD_vs_SPOD_L2_ntF1} and \ref{fig:MJO_POD_vs_SPOD_L2_ntF3} are obtained using all the 366 frequencies computed for the SPOD approach. We note that the dimensions of the latent space for POD ranged from 1 (mode) to 5 (modes), achieving a compression of $[1\text{ to }5] / 115680$. The one for SPOD ranged from $1 \times 366 = 366$ to $5 \times 366 = 1830$, resulting in a data compression of $[366\text{ to }1830] / 115680$, that is two orders of magnitude larger than the one obtained using POD. 

\subsubsection{SPOD emulation: Sensitivity to number of frequencies}\label{subsec:mjo_spod_freqs}
\begin{figure}[H]
\centering
\begin{subfigure}{.49\textwidth}\label{fig:MJO_proj}
\includegraphics[width=1.\textwidth]{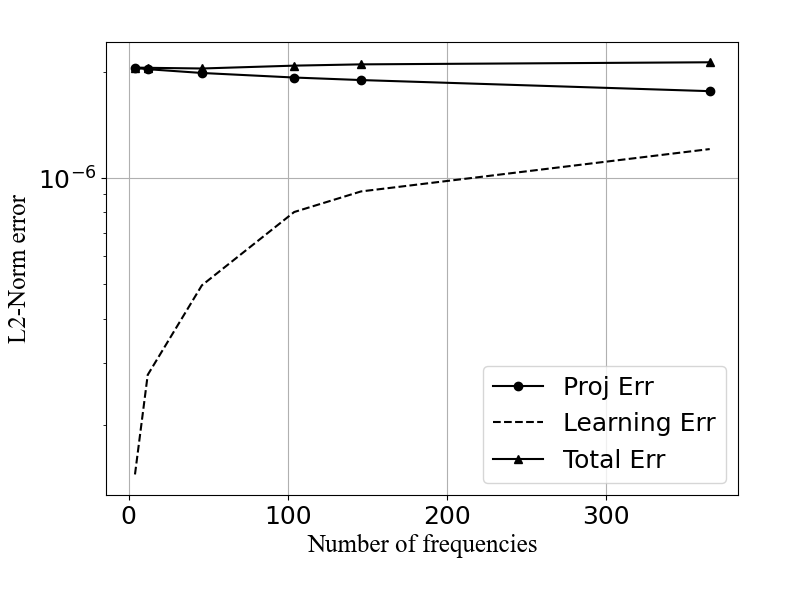}
\end{subfigure}
\begin{subfigure}{.49\textwidth}\label{fig:MJO_lstm}
\includegraphics[width=1.\textwidth]{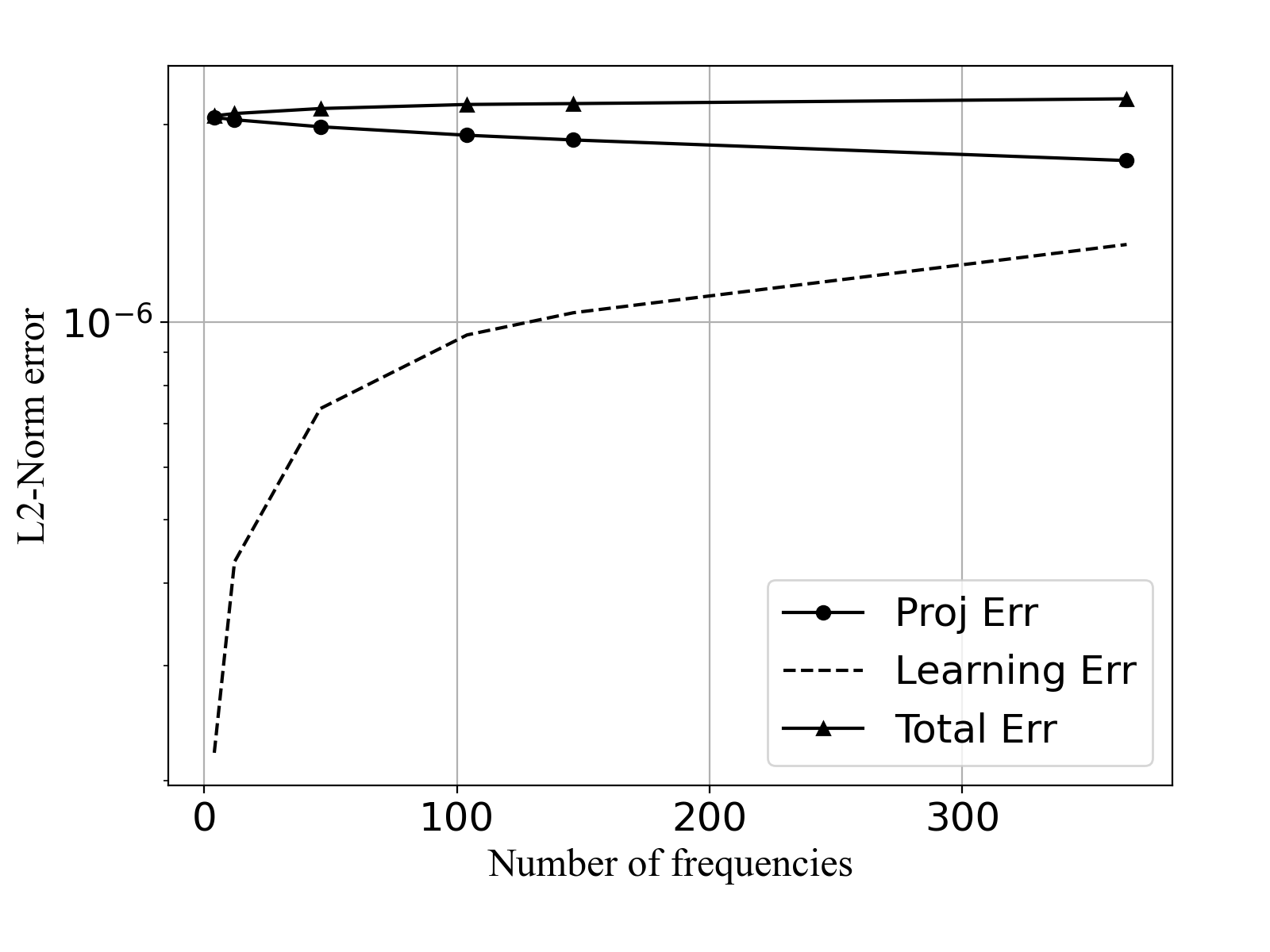}
\end{subfigure}
\caption{Projection, learning, and total $L_2$ errors plotted against the number of frequencies, retaining 3 modes, for forecast horizon $n_{\tau F} = 1$ (left plot), and $n_{\tau F} = 3$ (right plot).}
\label{fig:MJO_errorL2_freq}
\end{figure}

In this section, we show the sensitivity of the emulation strategy with respect to the number of frequencies retained for the SPOD reduction, similarly to what we showed in section~\ref{sec:comprJet}.
% Figure~\ref{fig:MJO_modes} shows the first mode for four distinct frequencies (periods). 
% %
% \begin{figure}[H]
% \centering
% \includegraphics[width=0.95\textwidth]{figures/MJO_modes.png}
% \caption{Visualizations of the first mode for four different periods: (a) $T=45.64$ , (b) $T = 40.57$, (c) $T=36.52$, and (d) $T=33.20$}
% \label{fig:MJO_modes}
% \end{figure}
% %
In detail, figure~\ref{fig:MJO_errorL2_freq} shows the L2 errors for different numbers of SPOD frequencies retained. As expected (and similar to the jet test case), a higher number of retained frequencies tends to reduce the projection error while, at the same time, affecting negatively the learning error. The two plots are referred to two different forecast horizons (12 and 36 hours, on the left and on the right, respectively).

\begin{table}[H]
\centering
{
\scriptsize
\begin{tabular}{?p{0.9cm}?p{1.1cm}|p{1.1cm}|p{1.1cm} ?p{1.1cm}|p{1.1cm}|p{1.1cm}?p{1.1cm}|p{1.1cm}|p{1.1cm}?}
 \hline
\textbf{$n_{\tau F} = 1$} &\multicolumn{3}{c?}{\textbf{Projection Error}} &\multicolumn{3}{c?}{\textbf{Learning Error}} &\multicolumn{3}{c?}{\textbf{Total Error}} \\
 \hline
 \hline
 Nr. of Freqs &  $L_{1}$ & $L_2$ & $L_{\infty}$  & $L_1$ & $L_2$ & $L_{\infty}$  & $L_1$ & $L_2$ & $L_{\infty}$ \\
 \hline
 4 & 3.23e-04  &  2.05e-06  &  3.18e-02  &  3.43e-05  &  1.46e-07  &  5.1e-04  &  3.23e-04  &  2.05e-06  &  3.19e-02 \\
 \hline
 12 & 3.25e-04  &  2.02e-06  &  3.16e-02  &  6.04e-05  &  2.77e-07   &  9.12e-04  &  3.26e-04  &  2.05e-06  &  3.19e-02 \\
 \hline
46 & 3.31e-04  &  1.96e-06  &  3.09e-02  &  1.41e-04  &   4.97e-07  &  1.92-03  &  3.41e-04  &  2.04e-06  &  3.18e-02 \\
\hline 
104 & 3.37e-04  &  1.91e-06  &  2.99e-02  &  2.01e-04  &   8.01e-07  &  2.46e-03  &  3.52-04  &  2.08e-06  &  3.18e-02  \\
\hline 
146 & 3.40e-04  &  1.85e-06  &  2.92e-02  &  2.25e-04  &  9.16e-07 &  3.50e-03  &  3.59-04  &  2.10e-06  &  3.18e-02 \\
 \hline
 365 & 3.41e-04  &  1.76e-06  &  2.68e-02  &  2.67e-04  &  1.21e-06  &  5.21e-03  &  3.68-04  &   2.13e-06  &  3.18e-02 \\
\hline 
\end{tabular}

\begin{tabular}{?p{0.9cm}?p{1.1cm}|p{1.1cm}|p{1.1cm} ?p{1.1cm}|p{1.1cm}|p{1.1cm}?p{1.1cm}|p{1.1cm}|p{1.1cm}?}
 \hline
\textbf{$n_{\tau F} = 3$} &\multicolumn{3}{c?}{\textbf{Projection Error}} &\multicolumn{3}{c?}{\textbf{Learning Error}} &\multicolumn{3}{c?}{\textbf{Total Error}} \\
 \hline
 \hline
 Nr. of Freqs &  $L_{1}$ & $L_2$ & $L_{\infty}$  & $L_1$ & $L_2$ & $L_{\infty}$  & $L_1$ & $L_2$ & $L_{\infty}$ \\
 \hline
 4 & 3.23e-04  &  2.05e-06  &  3.18e-02  &  5.19e-05  &  2.21e-07  &  7.53e-04  &  3.25e-04  &  2.06e-06  &  3.19e-02 \\
 \hline
 12 & 3.25e-04  &  2.03e-06  &  3.16e-02  &  1.04e-04  &  4.31e-07  &  1.42e-03  &  3.35e-04  &  2.08e-06  &  3.19e-02 \\
 \hline
46 & 3.31e-04  &  1.98e-06  &  3.09e-02  &  1.80e-04  &  7.38e-07  &  2.49e-03  &  3.56e-04  &  2.11e-06  &  3.18e-02 \\
\hline 
104 & 3.37e-04  &  1.92e-06  &  2.99e-02  &  2.32e-04  &  9.55e-07  &  3.65e-03  &  3.69e-04  &  2.14e-06  &  3.19e-02  \\
\hline 
146 & 3.40e-04  &  1.89e-06  &  2.92e-02  &  2.50e-04  &  1.03e-06  &  4.38e-03  &  3.72e-04  &  2.15e-06  &  3.19e-02 \\
 \hline
 365 & 3.41e-04  &  1.76e-06  &  2.68e-02  &  3.14e-04  &  1.31e-06  &  7.51e-03  &  3.88e-04  &  2.19e-06  &  3.19e-02 \\
 \hline 
\end{tabular}
}
\caption{$L_1$, $L_2$, and $L_{inf}$ norm errors evaluated retaining a different number of SPOD frequencies for $n_{\tau F} = 1$ (top) and $n_{\tau F}=3$.}
\label{tab:mjo_spod_errors_freqs}
\end{table}

\begin{figure}[H]
\centering
\begin{subfigure}{.49\textwidth}\label{fig:MJO_23-30_proj}
\includegraphics[width=\linewidth]{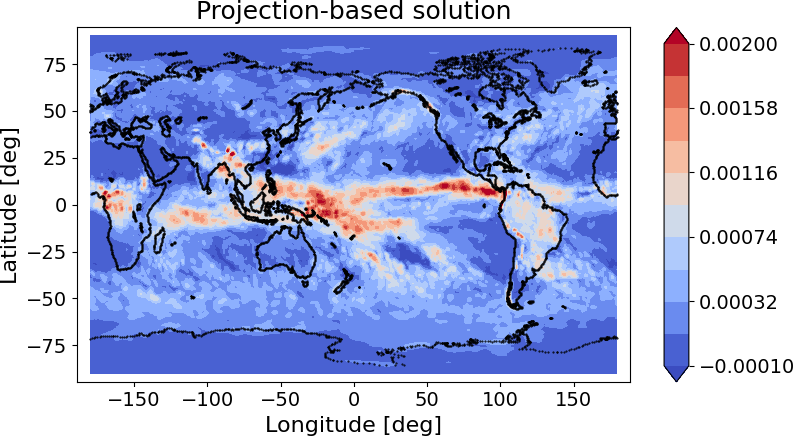}
\end{subfigure}
\begin{subfigure}{.49\textwidth}\label{fig:MJO_23-30_lstm}
\includegraphics[width=\linewidth]{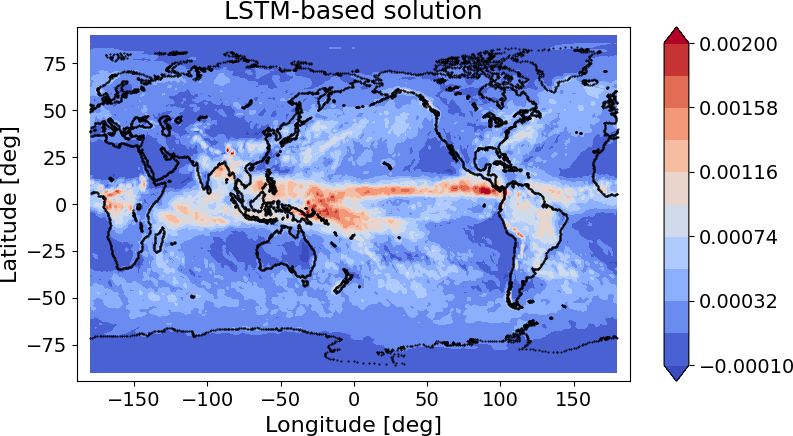}
\end{subfigure}

\vspace{0.3cm}
\begin{subfigure}{.49\textwidth}\label{fig:MJO_proj_30-50}
\includegraphics[width=\linewidth]{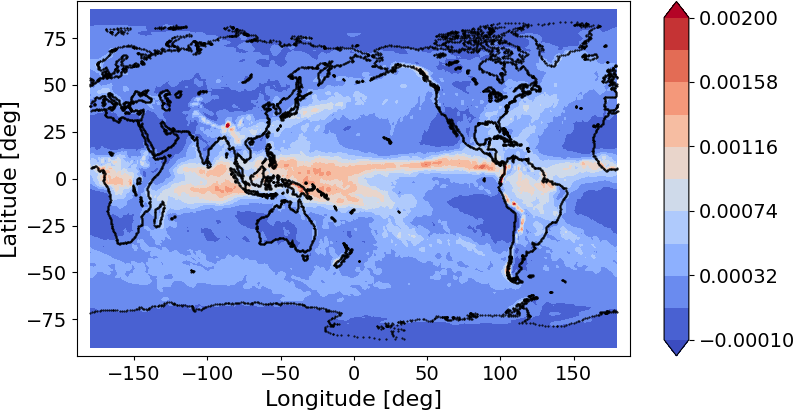}
\end{subfigure}
\begin{subfigure}{.49\textwidth}\label{fig:MJO_lstm_30-50}
\includegraphics[width=\linewidth]{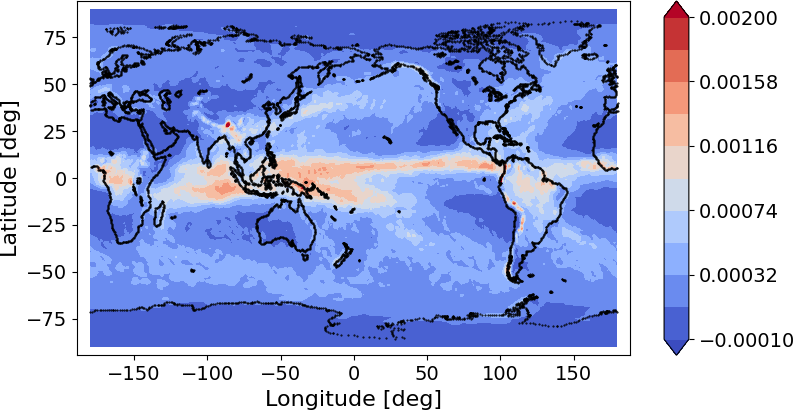}
\end{subfigure}

\vspace{0.3cm}
\begin{subfigure}{.49\textwidth}\label{fig:MJO_60-365_proj}
\includegraphics[width=\linewidth]{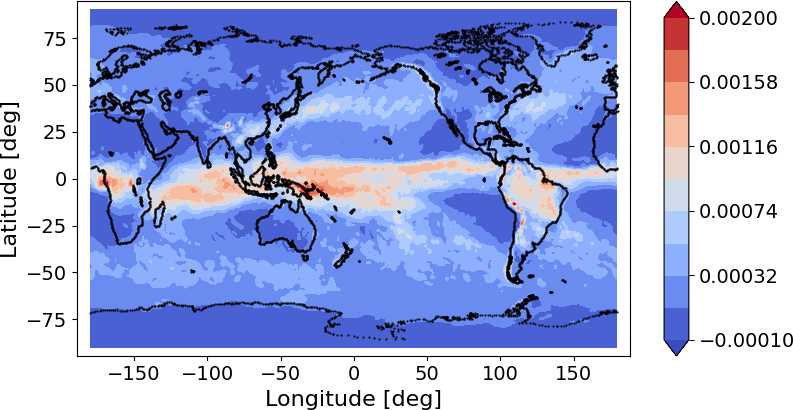}
\end{subfigure}
\begin{subfigure}{.49\textwidth}\label{fig:MJO_60-365_lstm}
\includegraphics[width=\linewidth]{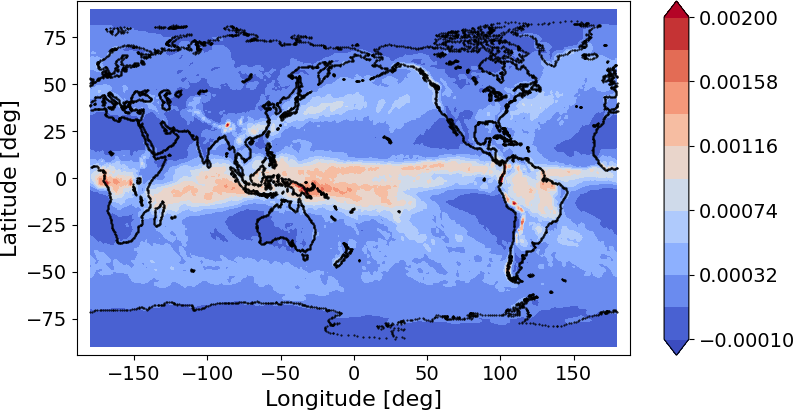}
\end{subfigure}
\caption{Comparison between the solution obtained with the coefficients projected onto the reduced space (left) and the one evaluated starting from the one learnt by the LSTM using $n_{\tau F} = 1$ (right) for three period ranges: 23 to 30 (top), 30 to 50 (middle), and 60 to 365 days. The results have been obtained using the same number of frequencies (i.e., 6), retaining 3 modes, and considering the same snapshot.}
\label{fig:MJO_SPOD_redFreqSnap}
\end{figure}

The prediction of the 300-th snapshot belonging to the testing set is reported in figure~\ref{fig:MJO_SPOD_redFreqSnap}. In particular, in all plots we retained the same number of frequencies (6), but we used different  period intervals: one for low frequencies (70--365 days), one for intraseasonal frequencies (30--50 days), and one for high frequencies (23--30 days). The left plots show the projection-based solutions, while the right plots show the LSTM predictions for forecast horizon $n_{\tau F} = 1$ (12 hours). The LSTM forecast, is qualitatively similar to the projection-based solution, as also highlighted by the relatively small L2 errors reported in figure~\ref{fig:MJO_errorL2_freq}.
%

%
% \begin{figure}[H]
% \centering
% \begin{subfigure}{.49\textwidth}\label{fig:MJO_23-30}
% \includegraphics[width=\linewidth]{figures/MJO_23-30.png}
% \caption{Period:23-30 days}
% \end{subfigure}
% \begin{subfigure}{.49\textwidth}\label{fig:MJO_70-365}
% \includegraphics[width=\linewidth]{figures/MJO_60-365.png}
% \caption{Period: 70-365 days}
% \end{subfigure}
% \caption{Projection solution (top) vs LSTM prediction (bottom) of the same time instant using two different reduced periods: 23--30 days (left) and 70--365 (right). Both results have been obtained using the same number of frequencies (i.e., 6), and retaining 3 modes.}
% \label{fig:SPOD_redFreq}
% \end{figure}
% \begin{figure}[H]
% \centering
% \includegraphics[width=1.1\linewidth]{figures/MJO_60-365_2.png}
% \caption{Projection solution (top) vs LSTM prediction (bottom) of the same time instant using two different reduced periods: 23--30 days (left) and 70--365 (right). Both results have been obtained using the same number of frequencies (i.e., 6), and retaining 3 modes.}
% \label{fig:SPOD_redFreq2}
% \end{figure}
%

%%%%%%%%%%%%%%%%%%%%%%%%%%%%%%%%%%%%%%%%%%%%%
%
\section{Conclusions}\label{sec:conclusions}
%
%%%%%%%%%%%%%%%%%%%%%%%%%%%%%%%%%%%%%%%%%%%%%
This work has been driven by recent advancement in SPOD-based reduced order modeling for quasi-stationary data. In particular, we presented a surrogate model for the SPOD latent space via recurrent neural networks based on three fundamental blocks: a) extraction of the SPOD modes from a database of previously computed snapshots, b) estimation of coefficients which are able to reconstruct the original data realization when applied to the SPOD modes through an oblique projection, and c) prediction of future evolution in time of the data with the aid of LSTM neural networks. The resulting model is non-intrusive, stable, and it is able to reproduce the dynamics of the underlying data.

Moreover an extensive comparison between a SPOD-based model and its POD-based counterpart has been carried out on two distinct test cases. The numerical results show that SPOD modes provide a larger and more refined basis and this contributes to reduce the projection error with respect to a standard approach which relies on POD. The main drawback is a consistent computational overhead with respect to POD since we consider the evolution in time of multiple dynamics characterized by small frequency ranges. The large size of the latent space affects also the quality of the final predicted solution: the tests we have presented show that the learning error in SPOD tends to grow with the number of modes since numerical accumulation of errors occurs. In addition plots of the evolution of the SPOD coefficients in time show a behavior which is less harmonic-like and regular than the one of POD. The learning error could be probably reduced with a refined tuning of the hyperparameters of the LSTM networks, but the general considerations reported above still hold; the refinement of the parameters was out of purpose of the present work, which is focused on the methodological aspects. 

The main advantage of SPOD over the POD model consists into its versatility: in fact the user is allowed to provide a range of frequency and follow the evolution of phenomena which are relevant for the selected time scale.
As it has been pointed out in previous work \cite{nekkanti2021frequency,chu2021stochastic}, this is fundamental when one has to carry out frequency analysis, data de-noising, or prewhitening , procedures which are common both in industry and in environmental sciences.
The dependence of the L2-norm errors for different numbers of retained frequencies have been studied for all the test cases proposed and, as expected, a higher number of frequencies is able to produce accurate reconstructions of the full data field but, as for the number of modes, a bigger latent space implies a bigger learning error.

The framework proposed is general and can be adopted in conjunction with alternative neural-network architectures that may improve the learning error. 
Extensions of the present work could aim at improving the efficiency and the accuracy of the model; this can be achieved either by investigating alternative networks and strategies for automatic optimization of the hyperparameters or investigating the application of sparse identification of nonlinear dynamics techniques to SPOD-based reduced models.
In addition, a comparison against nonlinear methods for seeking latent space with the aid of AI-based architecture (e.g. autoencoders or diffusion maps) will be carried out in upcoming works.
Furthermore, understanding the interplay between the projection, learning and total errors constitutes an additional and crucial area of research for non-intrusive surrogate models like the one presented in this paper.

\section*{Acknowledgements}
This work was partially supported by MIUR (Italian ministry for university and research) through FARE-X-AROMA-CFD project, P.I. Prof. Gianluigi Rozza. Gianmarco Mengaldo acknowledges support from NUS startup grant R-265-000-A36-133. This material is partially based upon work supported by the U.S. Department of Energy (DOE), Office of Science, Office of Advanced Scientific Computing Research, under Contract DE-AC02-06CH11357.
Prof. Gianluigi Rozza acknowledges the support by European Union Funding for Research and Innovation – Horizon 2020 Program – in the framework of European Research Council Executive Agency: Consolidator Grant H2020 ERC CoG 2015 AROMA-CFD project 681447 “Advanced Reduced Order Methods with Applications in Computational Fluid Dynamics”.

%% New version of the num-names style
\bibliographystyle{elsarticle-num-names}
\bibliography{references.bib}

\end{document}